\documentclass[11pt,a4paper]{amsart}

\usepackage{amsmath}
\usepackage{mathtools}
\usepackage{amsthm}
\usepackage{amssymb}
\usepackage{hyperref}

\newcommand*\fullref[3][\relax]{%
  \ifdefined\hyperref%
    {\hyperref[#3]{#2\penalty 200\ \ref*{#3}#1}}%
  \else%
    {#2\penalty 200\ \relax\ref{#3}#1}%
  \fi%
}

\newcommand*{\defterm}[1]{\emph{#1}}
\DeclarePairedDelimiter{\parens}{\lparen}{\rparen}
\DeclarePairedDelimiter{\set}{\{}{\}}
\DeclarePairedDelimiterX{\gset}[2]{\{}{\}}{\,#1:#2\,}
\newcommand*{\nset}{\mathbb{N}}
\newcommand*{\rset}{\mathbb{N}}
\newcommand*{\emptyword}{\varepsilon}

\DeclarePairedDelimiterX{\pres}[2]{\langle}{\rangle}{#1\,\delimsize\vert\,\mathopen{}#2}
\newcommand*{\drel}[1]{\mathcal{#1}}

\newcommand*{\hypo}{{\mathsf{hypo}}}
\newcommand*{\sylv}{{\mathsf{sylv}}}
\newcommand*{\sylvsharp}{{{\mathsf{sylv}}\smash{{}^{\mathsf{\#}}}}}
\newcommand*{\baxt}{{\mathsf{baxt}}}
\newcommand*{\stal}{{\mathsf{stal}}}
\newcommand*{\taig}{{\mathsf{taig}}}

\newcommand*{\lps}{{\mathsf{lPS}}}
\newcommand*{\rps}{{\mathsf{rPS}}}

\newcommand*{\hypocong}{\equiv_\hypo}
\newcommand*{\stalcong}{\equiv_\stal}
\newcommand*{\sylvcong}{\equiv_\sylv}
\newcommand*{\sylvsharpcong}{\equiv_\sylvsharp}
\newcommand*{\taigcong}{\equiv_\taig}
\newcommand*{\baxtcong}{\equiv_\baxt}

\newcommand*{\rpscong}{\equiv_\rps}

\newcommand*{\plr}{\mathrm{P}^\rightarrow}
\newcommand*{\prl}{\mathrm{P}^\leftarrow}
\newcommand*{\plit}{\mathrm{P}}

\newcommand*{\phypo}[2][]{\plr_{\hypo}\parens[#1]{#2}}

\newcommand*{\psylv}[2][]{\prl_{\sylv}\parens[#1]{#2}}

\newcommand*{\psylvsharp}[2][]{\plr_{\sylvsharp}\parens[#1]{#2}}

\newcommand*{\pbaxt}[2][]{\plit_{\baxt}\parens[#1]{#2}}

\newcommand*{\ptaig}[2][]{\prl_{\taig}\parens[#1]{#2}}

\newcommand*{\pstal}[2][]{\prl_{\stal}\parens[#1]{#2}}

\newcommand*{\prps}[2][]{\plr_{\rps}\parens[#1]{#2}}

\newcommand*{\evlit}{{\mathrm{ev}}}
\newcommand*{\ev}[2][]{\evlit\parens[#1]{#2}}

\newcommand*\var[1]{\mathbf{#1}}

\usepackage{tikz}
\usetikzlibrary{calc}

\usetikzlibrary{shapes.geometric} 
\usetikzlibrary{shapes.misc}  

\tikzset{
  bst/.style={
    standard/.style={
      font=\small,
      draw=gray,
      rounded rectangle,
      minimum width=4.5mm,
      minimum height=4.5mm,
      inner xsep=0mm,
      inner ysep=1mm,
      outer sep=0mm,
      line width=.5pt,
    },
    empty/.style={
      minimum width=3mm,
      minimum height=3mm,
    },
    triangle/.style={
      isosceles triangle,
      isosceles triangle apex angle=60,
      shape border rotate=90,
      rounded corners=2mm,
      minimum width=8mm,
      inner xsep=0mm,
      inner ysep=.5mm
    },
    blank/.style={
      draw=none,
    },
    nodecount/.style={
      blank,
      font=\scriptsize,
    },
    every node/.style={standard},
    every child/.style={draw=black,line width=.6pt},
    level distance=10mm,
    level 1/.style={sibling distance=60mm},
    level 2/.style={sibling distance=30mm},
    level 3/.style={sibling distance=15mm},
  },
  medbst/.style={
    bst,
    level distance=10mm,
    level 1/.style={sibling distance=15mm},
    level 2/.style={sibling distance=15mm},
    level 3/.style={sibling distance=15mm},
  },
  smallbst/.style={
    bst,
    level distance=8mm,
    level 1/.style={sibling distance=10mm},
    level 2/.style={sibling distance=10mm},
    level 3/.style={sibling distance=10mm},
  },
  tinybst/.style={
    bst,
    level distance=5mm,
    level 1/.style={sibling distance=8mm},
    level 2/.style={sibling distance=8mm},
    level 3/.style={sibling distance=8mm},
    every node/.append style={
      font=\footnotesize,
    },
    triangle/.append style={
      rounded corners=1mm,
      minimum width=7mm,
      inner xsep=-.5mm,
    },
  },
  microbst/.style={
    bst,
    standard/.append style={
      font=\scriptsize,
      minimum width=3mm,
      minimum height=3mm,
      inner ysep=.25mm,
    },
    level distance=3mm,
    level 1/.style={sibling distance=6mm},
    level 2/.style={sibling distance=6mm},
    level 3/.style={sibling distance=6mm},
  },
  nanobst/.style={
    bst,
    standard/.append style={
      font=\tiny,
      minimum width=2mm,
      minimum height=2mm,
      inner ysep=.25mm,
    },
    level distance=2mm,
    level 1/.style={sibling distance=4mm},
    level 2/.style={sibling distance=4mm},
    level 3/.style={sibling distance=4mm},
  },
}
\usetikzlibrary{matrix}

\tikzset{
  pretableaumatrix/.style={
    ampersand replacement=\&,
    matrix of math nodes,
    outer sep=1mm,
    inner sep=0mm,
    anchor=center,
    row sep={between borders,-\pgflinewidth},
    column sep={between borders,-\pgflinewidth},
    dottedentry/.style={densely dotted},
    dashedentry/.style={densely dashed},
    spaceentry/.style={draw=none,execute at begin node=\null},
  },
  pretableaunode/.style={
    font=\small,
    draw=gray,
    sharp corners,
    rectangle,
    anchor=base,
    text height=3.75mm,
    text depth=1.25mm,
    minimum height=5mm,
    minimum width=5mm,
    inner sep=0mm,
    outer sep=0mm,
    doublewidth/.style={minimum width=10mm},
    footnotesize/.style={font=\footnotesize},
    scriptsize/.style={font=\scriptsize},
  },
  tableaumatrix/.style={
    pretableaumatrix,
    every node/.append style={
      pretableaunode,
    },
  },
  medtableaumatrix/.style={
    pretableaumatrix,
    every node/.append style={
      pretableaunode,
      font=\footnotesize,
      text height=2.75mm,
      text depth=.75mm,
      minimum height=3.5mm,
      minimum width=3.5mm
    },
  },
  smalltableaumatrix/.style={
    pretableaumatrix,
    every node/.append style={
      pretableaunode,
      font=\scriptsize,
      text height=1.85mm,
      text depth=.15mm,
      minimum height=2.5mm,
      minimum width=2.5mm,
    },
  },
  tinytableaumatrix/.style={
    pretableaumatrix,
    every node/.append style={
      pretableaunode,
      font=\tiny,
      text height=1.25mm,
      text depth=.15mm,
      minimum height=1.75mm,
      minimum width=1.75mm
    },
  },
  tableau/.style={
    baseline=-1.25mm,
    every matrix/.style={tableaumatrix},
  },
  medtableau/.style={
    baseline=-1.25mm,
    every matrix/.style={medtableaumatrix},
  },
  smalltableau/.style={
    baseline=-1.25mm,
    every matrix/.style={smalltableaumatrix},
  },
  preshapetableaumatrix/.style={
    pretableaumatrix,
    execute at end cell={\strut},
    every node/.append style={
      draw=black,
      anchor=base,
      inner sep=0mm,
      outer sep=0mm,
    },
    shadedentry/.style={fill=gray},
    darkshadedentry/.style={fill=darkgray},
  },
  medshapetableaumatrix/.style={
    preshapetableaumatrix,
    every node/.append style={
      text height=2.75mm,
      text depth=.75mm,
      minimum height=3.5mm,
      minimum width=3.5mm
    },
  },
  shapetableaumatrix/.style={
    ampersand replacement=\&,
    matrix of math nodes,
    outer sep=0mm,
    inner sep=0mm,
    anchor=base,
    row sep={between borders,-\pgflinewidth},
    column sep={between borders,-\pgflinewidth},
    execute at begin cell={\strut},
    every node/.append style={draw,anchor=base,text height=1mm,text depth=.5mm,minimum size=1.5mm,inner sep=0mm,outer sep=0mm},
  },
  shapetableau/.style={
    every matrix/.style={shapetableaumatrix},
  },
  topalign/.style={
    every matrix/.append style={name=maintableau,anchor=maintableau-1-1.base},
    baseline,
  },
}

\newcommand*\tableau[2][]{\tikz[tableau,#1]\matrix{#2};}

\newcommand*{\trop}{\mathbb{T}}
\newcommand*{\uttrop}[1]{\mathrm{UT}_{#1}(\trop)}
\newcommand*{\utS}[1]{\mathrm{UT}_{#1}(S)}

\newcommand*{\matS}[1]{\mathrm{M}_{#1}(S)}
\newcommand*{\ut}[2]{\mathrm{UT}_{#1}(#2)}
\newcommand*{\mat}[2]{\mathrm{M}_{#1}(#2)}
\newcommand*\supp[1]{\mathrm{supp}(#1)}

\newtheorem{theorem}{Theorem}[section]
\newtheorem{lemma}[theorem]{Lemma}

\newtheorem{corollary}[theorem]{Corollary}
\newtheorem{proposition}[theorem]{Proposition}

\theoremstyle{definition}

\newtheorem{remark}[theorem]{Remark}

\newtheorem{algorithm}[theorem]{Algorithm}

\numberwithin{equation}{section}

\allowdisplaybreaks
\begin{document}

\title[Representations and Identities of Plactic-like Monoids]{Representations and Identities of Plactic-like Monoids}

\author[A. J. Cain]{Alan J. Cain}
\address[A. J. Cain]{%
  Centro de Matem\'{a}tica e Aplica\c{c}\~{o}es\\
  Faculdade de Ci\^{e}ncias e Tecnologia\\
  Universidade Nova de Lisboa\\
  2829--516 Caparica\\
  Portugal
}
\email{%
  a.cain@fct.unl.pt
}
\thanks{The first and fourth authors were partially supported by the Funda\c{c}\~{a}o para a Ci\^{e}ncia e a
  Tecnologia (Portuguese Foundation for Science and Technology) through the projects {\scshape UIDB}/{\sc MAT}/2020
  (Centro de Matem\'{a}tica e Aplica\c{c}\~{o}es) and {\scshape PTDC}/{\scshape MAT-PUR}/31174/2017.}

\author[M. Johnson]{Marianne Johnson}
\address[M. Johnson]{%
  Department of Mathematics\\
  University of Manchester\\
  Manchester M13 9PL\\
  UK
}
\email{Marianne.Johnson@manchester.ac.uk}

\author[M. Kambites]{Mark Kambites}
\address[M. Kambites]{%
  Department of Mathematics\\
  University of Manchester\\
  Manchester M13 9PL\\
  UK
}
\email{Mark.Kambites@manchester.ac.uk}

\author[A. Malheiro]{Ant\'onio Malheiro}
\address[A. Malheiro]{%
  Departamento de Matem\'{a}tica \&\ Centro de Matem\'{a}tica e Aplica\c{c}\~{o}es\\
  Faculdade de Ci\^{e}ncias e Tecnologia\\
  Universidade Nova de Lisboa\\
  2829--516 Caparica\\
  Portugal
}
\email{%
  ajm@fct.unl.pt
}

\begin{abstract}
  We exhibit faithful representations of the hypoplactic, stalactic, taiga, sylvester, Baxter and right patience sorting monoids of each finite
  rank as monoids of upper triangular matrices over any semiring from a large class including the tropical semiring and fields of characteristic
  $0$. By analysing the image of
  these representations, we show that the variety generated by a single hypoplactic (respectively, stalactic or taiga)
  monoid of rank at least $2$ coincides with the variety generated by the natural numbers together with a fixed finite
  monoid $\mathcal{H}$ (respectively, $\mathcal{F}$) forming a proper subvariety of the variety generated by the plactic
  monoid of rank $2$.
\end{abstract}


\maketitle

\section{Introduction}
\label{sec:intro}

The plactic monoid (the monoid of Young tableaux) is famous for its connections to such diverse areas as symmetric
functions \cite{macdonald_symmetric}, representation theory and algebraic combinatorics
\cite{fulton_young,lothaire_algebraic}, Kostka--Foulkes polynomials \cite{lascoux_plaxique,lascoux_foulkes}, and musical
theory \cite{jedrzejewski_plactic}. Its finite-rank versions were shown to have faithful tropical representations by the
second and third authors \cite[Theorem~2.8]{johnson_tropical}. An important consequence of these representations, which
are specifically representations using upper triangular tropical matrices, is that each finite-rank plactic monoid
satisfies a non-trivial semigroup identity \cite[Theorem~3.1]{johnson_tropical}, which had been an actively-studied
question \cite{kubat_identities}. The dimension of the representation and thus the lengths of the resulting identities
are dependent on the rank of the monoid. The first and fourth authors, together with Kubat, Klein, and Okni\'{n}ski,
showed that the rank-$n$ plactic monoid does not satisfy any non-trivial identity of length less than or equal to $n$,
which implies that there is no single non-trivial identity satisfied by all finite-rank plactic monoids, and, moreover,
that the infinite-rank plactic monoid does not satisfy any non-trivial semigroup identity
\cite{ckkmo_placticidentity}. Kubat and Okni\'{n}ski \cite{kubat_plactic}, and C\'{e}do and the same authors
  \cite{cedo_irreducible} also studied representations over a field of the plactic algebra of ranks $3$ and $4$ (that
  is, the monoid ring of the plactic monoid of ranks $3$ or $4$ over the same field).

Plactic monoids belong to a family of `plactic-like' monoids which are connected with combinatorics and whose elements can be
identified with combinatorial objects. Other in the family include the hypoplactic monoids, whose elements are quasi-ribbon
tableaux and which play a role in the theory of quasi-symmetric functions analogous to that of the plactic monoid for
symmetric functions \cite{krob_noncommutative4,krob_noncommutative5,novelli_hypoplactic}; the sylvester and
$\#$-sylvester monoids, whose elements are respectively right strict and left strict binary search trees
\cite{hivert_sylvester}; the taiga monoids, whose elements are binary search trees with multiplicities
\cite{priez_lattice}; the stalactic monoids, whose elements are stalactic tableaux
\cite{hivert_commutative,priez_lattice}; the Baxter monoids, whose elements are pairs of twin binary search trees
\cite{giraudo_baxter,giraudo_baxter2} and which are linked to the theory of Baxter permutations; and the left and right
patience sorting monoids \cite{rey_algebraic,cms_patience1}, whose elements are patience sorting tableaux.

For each species of plactic-like monoid, with the exception of the patience sorting monoids, there exist fixed
identities satisfied by all finite- and infinite-rank monoids of that species \cite{cm_identities}. For the left
patience sorting monoid and its finite-rank versions, only the rank-$1$ monoid satisfies a non-trivial identity (since
it is commutative); for higher ranks, it contains a free submonoid of rank $2$ \cite[Corollary~4.4]{cms_patience1}. For the
right patience sorting monoid, the situation is similar to the plactic monoid: the rank-$n$ monoid satisfies no identity
of length less than of equal to $n$, and consequently the infinite-rank monoid satisfies no identity, but there are
identities, dependent on rank, satisfied by the finite-rank right patience sorting monoids.

This state of knowledge naturally raises the question of whether these plactic-like monoids admit faithful tropical
representations. Left patience sorting monoids of rank greater than $2$ are immediately excluded by their free submonoid
of rank $2$: it is known that finitely generated semigroups of tropical matrices have polynomial growth
\cite{dAllesandro03} and therefore cannot contain free submonoids of rank $2$ or more. In this paper, we show that, with
this exception, faithful tropical representations exist for all the plactic-like monoids mentioned above, namely the
hypoplactic, sylvester, $\#$-sylvester, Baxter, stalactic, taiga, and right patience sorting monoids. In fact we show that each of these monoids can be faithfully represented by upper triangular
matrices over any unital semiring with zero containing an element of infinite multiplicative order.

The paper is structured as follows. In Section 2 we outline some preliminary material on words, semigroups,
representations and identities. Each of \fullref{Sections}{sec:hypo}--\ref{sec:rps} studies a family of plactic-like
monoids. Each monoid in a given family is associated to a class of combinatorial object, and arises from an algorithm
that inserts a symbol into such an object. Starting from an empty object, it is therefore possible to compute a
combinatorial object from a word, and the elements of the monoid are equivalence classes of words that correspond to the
same object. We shall show that these monoids admit faithful representations by upper triangular matrices over
certain semirings.  \fullref{Section}{sec:hypo} concerns quasi-ribbon tableaux and representations of the hypoplactic
monoids; \fullref{Section}{sec:stal} concerns stalactic tableaux and representations of the stalactic monoids;
\fullref{Section}{sec:taig} concerns binary search trees with multiplicities and representations of the taiga monoids;
\fullref{Section}{sec:sylv} concerns binary search trees and representations of the sylvester and Baxter monoids;
and \fullref{Section}{sec:rps} concerns patience sorting tableaux and representations of the right patience sorting monoids.

\section{Preliminaries}
\label{sec:prelim}

\subsection{Words} We write $\nset$ for the set of positive integers and $\nset_0$ for the set of non-negative
integers. For $n \in\nset$ we write $[n]$ to denote the set $\set{1, \ldots n}$. For $i,j \in \nset$ we write $[[i,j]]$
to denote the set $\gset{k \in \nset}{\min(i,j) \leq k \leq \max(i,j)}$, or simply $[i,j]$ in the case where $i<j$, and
refer to such subsets as \emph{intervals}.  For a non-empty subset $X \subseteq \nset$ we write $X^*$ to denote the free
monoid generated by the set $X$, that is the set of all words on the (possibly ordered) alphabet $X$, where $\emptyword$
denotes the empty word.  For $w \in \nset^*$ we write $|w|$ to denote the length of the word $w$, and for each
$i \in \nset$ we write $|w|_i$ to denote the number of occurrences of the letter $i$ in $w$.

Each word $w \in \nset^*$ determines a function $\nset \rightarrow \nset_0$ via $x \mapsto |w|_x$ called the
\defterm{content} or \defterm{evaluation} of $w$, denoted $\ev{w}$ and a subset
$\supp{w} = \gset{x \in \nset}{|w|_x \neq 0} \subset \nset$, called the \defterm{support} of $w$. All of the
plactic-like monoids considered in this paper have a canonical generating set with the property that all words representing a given element have the
same content (and thus the same support). Hence it makes sense to define, for an element $m$ of the monoid, $\ev{m}$ and
$\supp{m}$ to be the content and support of the words representing $m$ with respect to the canonical generators.

\subsection{Matrix representations over semirings}
\label{subsec:trop-rep}

Throughout this paper $S$ will be a commutative unital semiring with zero denoted by $0_S$, unit denoted by $1_S$,
containing an element of infinite multiplicative order. Of particular interest is the \defterm{tropical semiring} $\trop$, which is the set
$\rset \cup\set{-\infty}$ under the operations $a \oplus b = {\rm max}(a,b)$ and $a \otimes b = a+b$, where we define
$-\infty + a=a+ -\infty = -\infty$. Notice that $0_{\trop}=-\infty$, $1_\trop = 0$ and all other elements have infinite
multiplicative order.

We write $\mat{n}{S}$ to denote the monoid of all $n \times n$ matrices with entries from $S$ under the matrix
multiplication induced from operations of $S$ in the obvious way. The $n \times n$ \textit{identity matrix} (with all diagonal equal to $1_S$ and all other entries equal to $0_S$) and \textit{zero matrix} (with all entries equal to $0_S$) are respectively the identity element and zero element in $\mat{n}{S}$. We say that $A \in \mat{n}{S}$ is \emph{upper triangular}
if $A_{i,j}=0_S$ for all $i>j$, and write $\ut{n}{S}$ for the submonoid of $n \times n$ upper triangular matrices over
$S$. If $T$ is a finite set we write $\mat{T}{S}$ for the semigroup of matrices with rows and columns indexed by
elements of $T$; this is of course isomorphic to $\mat{|T|}{S}$ but it  is often convenient to index entries
by elements of a particular finite set.

From now on fix an element of infinite multiplicative order $\alpha \in S$ and let
  $c_n : [n]^* \rightarrow \ut{n}{S}$ be the homomorphism extending the map defined for $x \in [n]$ by
\[
  c_n(x)_{p,q} =
  \begin{cases}
    \alpha       & \text{if $p=q=x$,}    \\
    1_S       & \text{if $p=q\neq x$,} \\
    0_S & \text{otherwise.}      \\
  \end{cases}
\]%
The image of this morphism is the (commutative) semigroup of diagonal matrices, with entries from
$\gset{\alpha^i}{i \in \nset_0}$ on the diagonal. Since $\alpha$ is an element of infinite multiplicative order,
this image is isomorphic to $n$ copies of the monoid $(\nset_0, +)$.  Two words $w, v \in [n]^*$ have the same
image under $c_n$ if and only if they have the same content. (Of course, if $S$ has the stronger property that its
multiplicative monoid contains a free commutative monoid of each finite rank $n$ (as is the case for the tropical
semiring, for example), then we can instead construct a $1$-dimensional representation of $[n]^*$ which records the
content of a word.)

\subsection{Identities}

A \defterm{semigroup identity} is a formal equality $u = v$ where $u$ and $v$ are non-empty words over some alphabet of variables
$X$. An identity is \defterm{non-trivial} if $u$ and $v$ are not equal as words. Such a semigroup identity $u = v$ is
\defterm{satisfied} by a semigroup $S$ if, for every homomorphism $\phi : X^+ \to S$, the equality $\phi(u) = \phi(v)$
holds in $S$. A semigroup identity is \defterm{balanced} (or \defterm{multihomogeneous}) if $\ev{u} = \ev{v}$.

It is easy to see that any identity satisfied by a semigroup containing a free submonoid of rank $1$ must be balanced:
to see that $u$ and $v$ must contain the same number of the variable $x$, consider the homomorphism sending $x$ to the
generator of the free submonoid and all other variables to the identity element of the monoid. All of the plactic-like
monoids considered in this paper contain free submonoids of rank $1$.

On the other hand, the monoid variety defined by all balanced identities is $\var{Comm}$, the class of commutative
monoids.

\section{Hypoplactic monoid}
\label{sec:hypo}

\subsection{The hypoplactic monoid}

A \defterm{quasi-ribbon tableau} is a planar diagram consisting of a finite array of adjacent symbols from $\nset$,
with each symbol lying either to the east or the south of the previous symbol, and with the property that symbols lying
in the same `row' form a non-decreasing sequence when read from left to right and symbols lying in the same `column' are
strictly increasing when read from top to bottom. An example of a quasi-ribbon tableau is:
\begin{equation}
  \label{eq:qrteg}
  \tableau
  {
    1 \& 1 \& 2 \&   \&   \&   \\
    \&   \& 3 \& 4 \& 4 \&   \\
    \&   \&   \&   \& 5 \&   \\
    \&   \&   \&   \& 6 \& 6 \\
  }.
\end{equation}
Notice that the same symbol cannot appear in two different rows of a quasi-ribbon tableau.

The insertion algorithm is as follows:

\begin{algorithm}[{\cite[\S~7.2]{krob_noncommutative4}}]
\label{alg:hypoinsertone}
~\par\nobreak
\textit{Input:} A quasi-ribbon tableau $T$ and a symbol $a \in \nset$.

\textit{Output:} A quasi-ribbon tableau $T\leftarrow a$.

\textit{Method:} If there is no entry in $T$ that is less than or equal to $a$, output the quasi-ribbon tableau obtained
by creating a new entry $a$ and attaching (by its top-left-most entry) the quasi-ribbon tableau $T$ to the bottom of $a$.

If there is no entry in $T$ that is greater than $a$, output the quasi-ribbon tableau obtained by creating a new entry $a$ and attaching
(by its bottom-right-most entry) the quasi-ribbon tableau $T$ to the left of $a$.

Otherwise, let $x$ and $z$ be the adjacent entries of the quasi-ribbon tableau $T$ such that $x \leq a < z$.
(Equivalently, let $x$ be the right-most and bottom-most entry of $T$ that is less than or equal to $a$, and let $z$ be
the left-most and top-most entry that is greater than $a$. Note that $x$ and $z$ could be either horizontally or
vertically adjacent.) Take the part of $T$ from the top left down to and including $x$, put a new entry $a$ to
the right of $x$ and attach the remaining part of $T$ (from $z$ onwards to the bottom right) to the bottom of the new
entry $a$, as illustrated here:
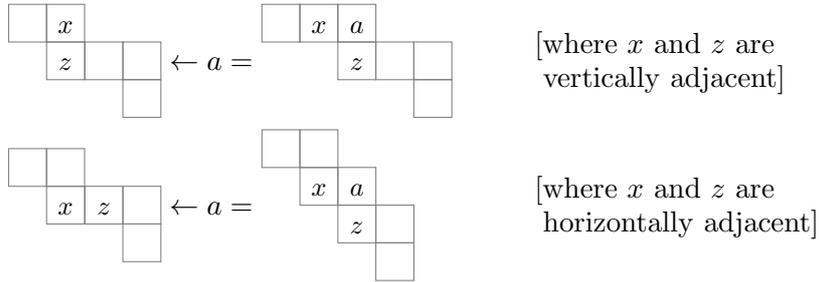
\begin{align*}
\tikz[tableau]\matrix
{
 \null \& x \&   \&     \\
 \& z \& \null \& \null \\
 \&   \&   \& \null     \\
}; \leftarrow a & =
\tikz[tableau]\matrix
{
 \null \& x \& a  \&          \\
 \& \& z \& \null \& \null    \\
 \& \& \&   \& \null          \\
}; &  &
\begin{array}{l}
\text{[where $x$ and $z$ are} \\
\text{\phantom{[}vertically adjacent]}
\end{array}
\displaybreak[0]              \\
\tikz[tableau]\matrix
{
 \null \& \null \&   \& \\
 \& x \& z \& \null \\
 \&   \&   \& \null \\
}; \leftarrow a &=
\tikz[tableau]\matrix
{
 \null \& \null \\
 \& x \& a \\
 \&  \& z \& \null \\
 \& \& \& \null \\
}; &&
\begin{array}{l}
\text{[where $x$ and $z$ are} \\
\text{\phantom{[}horizontally adjacent]}
\end{array}
\end{align*}
Output the resulting quasi-ribbon tableau.
\end{algorithm}

Thus one can compute, for any word $u \in \nset^*$, a quasi-ribbon tableau $\phypo{u}$ by starting with an empty
quasi-ribbon tableau and successively inserting the symbols of $u$, \emph{proceeding left-to-right} through the
word. Define the relation $\hypocong$ by
\[
u \hypocong v \iff \phypo{u} = \phypo{v}
\]
for all $u,v \in \nset^*$. The relation $\hypocong$ is a congruence, and the \defterm{hypoplactic monoid}, denoted
$\hypo$, is the factor monoid $\nset^*\!/{\hypocong}$; the \defterm{hypoplactic monoid of rank $n$}, denoted $\hypo_n$,
is the factor monoid $[n]^*/{\hypocong}$ (with the natural restriction of $\hypocong$). Each element $[u]_{\hypocong}$ (where $u \in \nset^*$) can be identified with the quasi-ribbon tableau $\phypo{u}$.

The monoid $\hypo$ is presented by
$\pres{\nset}{\drel{R}_\hypo}$, where
\begin{align*}
\drel{R}_\hypo ={} & \gset[\big]{(acb,cab)}{a,b,c \in \nset,\; a \leq b < c}               \\
                   & \cup \gset[\big]{(bac, bca)}{a,b,c \in \nset,\; a < b \leq c}         \\
                   & \cup \gset[\big]{(cadb,acbd)}{a,b,c \in \nset,\; a \leq b < c \leq d} \\
                   & \cup \gset[\big]{(bdac,dbca)}{a,b,c \in \nset,\; a < b \leq c < d};
\end{align*}
see \cite[\S~4.1]{novelli_hypoplactic} or \cite[\S~4.8]{krob_noncommutative4}. The monoid $\hypo_n$ is presented by
$\pres{[n]}{\drel{R}_\hypo}$, with the natural restrictions on the set of defining relations $\drel{R}_\hypo$.

\subsection{Constructing a faithful representation of the hypoplactic monoid}

For $u \in [n]$ and $1 \leq i <j \leq n$, let $H_{i,j}(u)$ denote the statement `$u$ contains $i$ and $j$, no symbol $k$
with $i < k < j$, and no scattered subword $ji$'. The following characterization of the hypoplactic monoid is a
consequence of \cite[Theorem 4.18 and Note 4.10]{novelli_hypoplactic}.

\begin{proposition}
  \label{prop:hypo_data}
  Let $n$ be a fixed positive integer and $u, v \in [n]^*$. The quasi-ribbon tableaux $\phypo{u}$ and $\phypo{v}$ are equal if and only if:
  \begin{enumerate}
  \item $u$ and $v$ have the same content; and
  \item for $1 \leq i<j \leq n$, $H_{i,j}(u) \iff H_{i,j}(v)$.
  \end{enumerate}
\end{proposition}

Since two words in the same hypoplactic class of rank $n$ have the same evaluation, they will have the same image under
the map $c_n$ from \fullref{Subsection}{subsec:trop-rep}. Therefore, by a slight abuse of notation, consider $c_n$ to be
a homomorphism from $\hypo_n$ to $\ut{n}{S}$, with image isomorphic to $n$ copies of the natural numbers.

Let
\[
  I =
  \begin{bmatrix}
    1_S & 1_S \\
    0_S & 0_S
  \end{bmatrix},
  \quad
  J =
  \begin{bmatrix}
    0_S & 0_S\\
    0_S & 1_S
  \end{bmatrix}.
\]
Let $K = JI$, $L = IJ$, and let $E$ be the $2\times2$ identity matrix. Note that $K$ is the $2 \times 2$ zero matrix. It
is easy to see that $\mathcal{H} = \set{E,I,J,K,L}$ is a submonoid of $\ut{2}{S}$, with presentation
$\pres{I, J}{I^2=I, J^2=J, IJI=JI=JIJ}$.  It can also be verified that $\mathcal{H}$ is isomorphic to the monoid
$\mathcal{C}_3$ of order-preserving and extensive transformations of the $3$-element chain, as studied in
\cite{Volkov04}.

For all $i,j \in \nset$ with $i<j$, consider the monoid homomorphism $h_{ij}:\nset^* \to \ut{2}{S}$ defined by
$i \mapsto I$, $j \mapsto J$, each $k$ with $i < k < j$ maps to $K$, and all other letters map to $E$. Note that the
image of $h_{ij}$ is $\mathcal{H}$. Straightforward calculation shows that for $w \in \nset^*$,
\begin{equation}
  \label{eq:hij-cases}
  h_{ij}(w) =
  \begin{cases}
    E & \text{if $w$ contains no symbols in the interval $[i,j]$;} \\
    I & \text{if $w$ contains $i$ and no other symbols from $[i,j$];} \\
    J & \text{if $w$ contains $j$ and no other symbols from $[i,j$];} \\
    L & \text{\parbox[t]{8cm}{\raggedright if $w$ contains $i$ and $j$, no other symbols from $[i,j]$, and no scattered subword $ji$; and}} \\
    K & \text{otherwise.}
  \end{cases}
\end{equation}

\begin{lemma}
  \label{lem:hypohoms}
  Let $n\geq 2$ be a fixed positive integer and let $1 \leq i < j \leq n$. The homomorphism
  $h_{ij}: \nset^* \rightarrow \ut{2}{S}$ factors to give a homomorphism from the hypoplactic monoid of rank $n$ to the
  monoid $\mathcal{H}$.
\end{lemma}

\begin{proof}
  Suppose that $u,v \in [n]^*$ are in the same hypoplactic class. Since $u$ and $v$ have the same content, it is
  immediate that either both or neither $u$ and $v$ have image $X$ where $X \in \set{E, I, J}$. Moreover, by
  \fullref{Proposition}{prop:hypo_data} $u$ has image $L$ if and only if $v$ has image $L$. Since the only other
  possible image is $K$, the result now follows.
\end{proof}

\begin{theorem}
  \label{thm_hypo}
  The hypoplactic monoid monoid $\hypo_n$ (respectively, $\hypo$) embeds into a direct product of $n$ copies
  (respectively, countably infinite copies) of $(\nset, +)$ with $\binom{n}{2}$ copies (respectively, countably infinite
  copies) of the finite monoid $\mathcal{H}$.
\end{theorem}

\begin{proof}
  The direct product of the morphisms $c_i$ for $i \in [n]$ (respectively, $i \in \nset$) and $h_{ij}$ for
  $i , j \in [n]$ (respectively, $i , j \in \nset)$ with $i < j$ give a morphism to the required monoid. The fact that
  this is an embedding now follows from \fullref{Proposition}{prop:hypo_data} and from observing from
  \eqref{eq:hij-cases} that the image of an element $w$ under $h_{ij}$ is equal to $L$ if and only if $w$ contains $i$
  and $j$, no other symbols from $[i,j]$, and no scattered subword $ji$.
\end{proof}

\begin{theorem}
  \label{thm_hypo_tropical}
  Let $S$ be a commutative unital semiring with zero containing an element of infinite multiplicative order.  The
  hypoplactic monoid of rank $n$ admits a faithful representation by upper triangular matrices of size $n^2$ over $S$
  having block-diagonal structure with largest block of size $2$ (or size $1$ if $n=1$).
\end{theorem}

\begin{proof}
  Since $(\nset, +)$ embeds in $\ut{1}{S}$ and $\mathcal{H}$ by definition embeds in $\ut{2}{S}$, this is immediate from
  \fullref{Theorem}{thm_hypo}.
\end{proof}

\subsection{Variety generated by the hypoplactic monoid}

The identities satisfied by the hypoplactic monoid have been completely characterized by the first and fourth authors and
Ribeiro \cite[Theorem~4.1]{cmr_idhypo}, but we deduce here some more information about the corresponding variety.

Let $J_k$ denote the set of identities $u=v$ with the property that $u$ and $v$ admit the same set of scattered subwords
of length at most $k$. Let $\var{Comm}$ denote the variety of commutative monoids (which is the variety with eqautional
theory given by the balanced identities), $\var{B}$ denote the variety of monoids generated by the bicyclic monoid, and
$\var{J}_k$ the variety of monoids determined by the set $J_k$.

The monoid varieties $\var{J}_k$ have been studied extensively \cite{straubing, pin, blanchet-sadri93, blanchet-sadri94,
  Volkov04}. By a result of Volkov \cite[Theorem~2]{Volkov04}, $\var{J}_k$ is generated by any one of: the monoid of
unitriangular Boolean matrices of rank $k+1$; the monoid reflexive binary relations on a set of size $k+1$; and most
importantly for our purposes, the monoid of order-preserving and extensive transformations of a chain with $k+1$
elements. In particular, the latter means that the five-element monoid $\mathcal{H}$ generates the variety
$\var{J}_2$. By a result of Tischenko \cite{tischenko} any five element monoid generates a finitely based variety, and
so certainly $\var{J}_2$ is finitely based. More generally, Blanchet-Sadri \cite{blanchet-sadri93, blanchet-sadri94} has
shown that the monoid variety $\var{J}_k$ is finitely based if and only if $k \leq 3$, providing a basis of identities
in those cases: a finite basis of identities for $\var{J}_2$ is $xyxzx = xyzx$, and $xyxy=yxyx$.

\begin{corollary}
  \label{corol:hypovar}
  Let $n\geq 2 $ be a fixed positive integer. The variety of monoids generated  by the hypoplactic monoid of rank $n$ is:
  \begin{enumerate}
  \item a proper subvariety of $\var{B}$;
  \item the join of $\var{Comm}$ and $\var{J}_2$;
  \item equal to the variety generated by the (infinite-rank) hypoplactic monoid.
  \end{enumerate}
\end{corollary}

\begin{proof}
  \begin{enumerate}
  \item We have seen that the hypoplactic monoid of rank $n$ embeds in the direct product of $n$ copies of
    $(\nset_0, +)$ (which embeds in $\uttrop{1}$) and $\binom{n}{2}$ copies of $\uttrop{2}$. Thus $\hypo_n$ is contained
    in the variety generated by $\uttrop{2}$. By \cite{daviaud_identities}, the latter is equal to the variety generated
    by the bicyclic monoid. That these varieties are distinct follows from the fact that the shortest identity satisfied
    by $\uttrop{2}$ has length $10$, whilst $\hypo_n$ satisfies the identity $xyxy=yxyx$
    \cite[Proposition~12]{cm_identities}.

  \item We begin by showing that the identities satisfied by the hypoplactic monoid of rank $n$ are precisely the
    balanced identities satisfied by the monoid $\mathcal{H}$. By \fullref{Theorem}{thm_hypo}, it is clear that the
    hypoplactic monoid of rank $n$ embeds in a direct product of copies of $(\nset_0, +)$ with copies of
    $\mathcal{H}$. Thus $\hypo_n$ satisfies every identity satisfied by both $(\nset_0, +)$ and $\mathcal{H}$. The
    identities satisfied by $(\nset_0, +)$ are precisely identities of the form $u=v$ where $u$ and $v$ have the same
    content, thus $\hypo_n$ satisfies all balanced identities satisfied by $\mathcal{H}$. On the other hand,
    $\mathcal{H}$ is an image of $\hypo_n$ under any of the homomorphisms $h_{ij}$, and so $\mathcal{H}$ satisfies every
    identity satisfied by $\hypo_n$. Finally, note that all identities satisfied by $\hypo_n$ are balanced.

    It now follows from the fact that $\mathcal{H}$ generates the variety $\var{J}_2$ \cite[Theorem~2]{Volkov04} that the variety generated by $\hypo_n$ is the join of
    $\var{Comm}$ and $\var{J}_2$.

  \item It is clear that each $\hypo_n$ embeds in $\hypo$, while it is known \cite[Proposition~3.6]{cmr_idhypo} that
    $\hypo$ can be embedded into a direct product of copies of $\hypo_2$. Thus, the varieties generated by all of these
    monoids coincide. \qedhere
  \end{enumerate}
\end{proof}

\section{Stalactic monoid}
\label{sec:stal}
\subsection{The stalactic monoid}

A \defterm{stalactic tableau} is a finite array of symbols from $\nset$ in which columns are top-aligned, and two symbols
appear in the same column if and only if they are equal. For example,
\begin{equation}
\label{eq:egstaltab1}
\tableau{
3 \& 1 \& 2 \& 6 \& 5 \\
3 \& 1 \&   \& 6 \& 5 \\
  \& 1 \&   \&   \& 5 \\
  \& 1 \\
}
\end{equation}
is a stalactic tableau. The insertion algorithm is very straightforward:

\begin{algorithm}
\label{alg:stalinsertone}
~\par\nobreak
\textit{Input:} A stalactic tableau $T$ and a symbol $a \in \nset$.

\textit{Output:} A stalactic tableau $T \leftarrow a$.

\textit{Method:} If $a$ does not appear in $T$, add $a$ to the left of the top row of $T$. If $a$ does appear in $T$,
add $a$ to the bottom of the (by definition, unique) column in which $a$ appears. Output the new tableau.
\end{algorithm}

Thus one can compute, for any word $u \in \nset^*$, a stalactic tableau $\pstal{u}$ by starting with an empty stalactic
tableau and successively inserting the symbols of $u$, \emph{proceeding right-to-left} through the word. For example
$\pstal{36113\allowbreak 5112565}$ is \eqref{eq:egstaltab1}. Notice that the order in which the symbols appear along the first row
in $\pstal{u}$ is the same as the order of the rightmost instances of the symbols that appear in $u$. Define the relation
$\stalcong$ by
\[
  u \stalcong v \iff \pstal{u} = \pstal{v}
\]
for all $u,v \in \nset^*$. The relation $\stalcong$ is a congruence, and the \defterm{stalactic monoid}, denoted $\stal$,
is the factor monoid $\nset^*\!/{\stalcong}$; The \defterm{stalactic monoid of rank $n$}, denoted $\stal_n$, is the factor
monoid $[n]^*/{\stalcong}$ (with the natural restriction of $\stalcong$). Each element $[u]_{\stalcong}$ (where
$u \in \nset^*$) can be identified with the stalactic tableau $\pstal{u}$. Note that if $T$ is a stalactic tableau
consisting of a single row (that is, with all columns having height $1$), then there is a unique word $u \in \nset^*$,
formed by reading the entries of $T$ left-to-right, such that $\pstal{u} = T$.

The monoid $\stal$ is presented by
$\pres{\nset}{\drel{R}_\stal}$, where
\[
\drel{R}_\stal = \gset[\big]{(bavb,abvb)}{a,b\in \nset,\; v \in \nset^*};
\]
see \cite[Example 3]{priez_lattice}. The monoid $\stal_n$ is presented by $\pres{[n]}{\drel{R}_\stal}$, where the set of
defining relations $\drel{R}_\stal$ is naturally restricted to $[n]^*\times [n]^*$.

\subsection{Constructing a faithful representation of the stalactic monoid}

For $u \in [n]$ and $i,j \in [n]$ with $i \neq j$, let $S_{i,j}(u)$ denote the statement `$u$ factors as $u=u'iu''$ with
where $u''$ contains $j$ but not $i$'. Notice that if $i$ and $j$ are in the support of $u$, then exactly one of
$S_{i,j}(u)$ and $S_{j,i}(u)$ is true. Otherwise both statements are false.

\begin{proposition}
  \label{prop:staldata}
  Let $n$ be a fixed positive integer and $u,v \in [n]^*$. The stalactic tableaux $\pstal{u}$ and $\pstal{v}$ are equal if and only if:
  \begin{enumerate}
  \item $u$ and $v$ have the same content; and
  \item for $1 \leq i < j \leq n$, $S_{i,j}(u) \iff S_{i,j}(v)$.
  \end{enumerate}
\end{proposition}

\begin{proof}
  It follows from the insertion algorithm that $\pstal{u}=\pstal{v}$ if and only if $u$ and $v$ have the same content
  and the order of the rightmost instances of the symbols that appear in $u$ is equal to the order of the rightmost
  instances of the symbols that appear in $v$. For each pair of elements $i,j$ in the support of $u$, the statement
  $S_{i,j}(u)$ (if $i<j$) or $S_{j,i}(u)$ (if $j<i$) can be used to determine whether $j$ occurs to the right of the
  right-most $i$, and hence (1) and (2) together imply $\pstal{u}=\pstal{v}$.

  Conversely, suppose that $\pstal{u}=\pstal{v}$. We immediately have (1) (and hence in particular, the support of $u$
  is equal to the support of $v$) and so it remains to show that (2) holds. For $i \in \supp{u}$ it follows from the
  insertion algorithm that if $u=u'iu''$ where $i \not\in \supp{u''}$, then $i$ occurs in the $(|\supp{u''}|+1)$-th
  column of the tableau (recalling here that the insertion algorithm reads the word from right to left) and the support
  of $u''$ is the set of entries in the first row preceding the symbol $i$. Thus $S_{i,j}(u)$ holds if and only if in
  the first row of $\pstal{u}=\pstal{v}$ symbol $j$ precedes symbol $i$ if and only if $S_{i,j}(v)$ holds.
\end{proof}

Two words in the same stalactic class of rank $n$ will have the same image under the content map $c_n$ from
\fullref{Subsection}{subsec:trop-rep}. Therefore, by a slight abuse of notation, consider $c_n$ as a homomorphism from
$\stal_n$ to $\ut{n}{S}$.

Let
\[
  I = \begin{bmatrix} 1_S&1_S\\ 0_S&0_S\end{bmatrix},\quad
  J = \begin{bmatrix} 1_S&0_S\\ 0_S&0_S\end{bmatrix},
\]
and let $E$ be the $2 \times 2$ identity matrix. Then $\mathcal{F} = \set{E,I,J}$ is isomorphic to the `flip-flop
monoid' (a two element right zero semigroup with an identity adjoined), presented by $\pres{I, J}{JI=I=I^2, IJ=J=J^2}$.

For all $i,j \in \nset$ with $i < j$, consider the monoid homomorphism $s_{ij} : \nset^* \to \ut{2}{S}$ defined by
$i \mapsto I$, $j \mapsto J$, and all other letters map to $E$. Note that the image of $s_{ij}$ is
$\mathcal{F}$. It is easy to see that for $w \in \nset^*$
\begin{equation}
  \label{eq:sij-cases}
  s_{ij} = \begin{cases}
    E & \text{if $i,j \notin \supp{w}$} \\
    I & \text{if $i \in \supp{w}$ and $w$ factors as $w = w'iw''$ with $j \notin \supp{w''}$} \\
    J & \text{if $j \in \supp{w}$ and $w$ factors as $w = w'jw''$ with $i \notin \supp{w''}$}
  \end{cases}
\end{equation}

\begin{lemma}
  \label{lem:stalhoms}
  Let $n\geq 2$ be a fixed positive integer and let $1 \leq i < j \leq n$. The homomorphism
  $s_{ij}: \nset^* \rightarrow \ut{2}{S}$ factors to give a homomorphism from the stalactic monoid of rank $n$ to the
  monoid $\mathcal{F}$.
\end{lemma}

\begin{proof}
  Let $u,v \in \nset^*$ be in the same stalactic class. Since $u$ and $v$ have the same content, it is immediate that
  either both or neither have image $E$. Moreover, by \fullref{Proposition}{prop:staldata}, $u$ has image $I$ if and
  only if $v$ has image $I$. Since $J$ is the only other possible image, this completes the proof.
\end{proof}

\begin{theorem}
  \label{thm_stal}
  The stalactic monoid $\stal_n$ (respectively, $\stal$) embeds into a direct product of $n$ copies (respectively,
  countably infinite copies) of $(\nset_0, +)$ and $\binom{n}{2}$ copies (respectively, countably infinite copies) of
  the finite monoid $\mathcal{F}$.
\end{theorem}

\begin{proof}
  The direct product of the morphisms $c_i$ for $i \in [n]$ (respectively, $i \in \nset$) and $s_{ij}$ for
  $i, j \in [n]$ (respectively, $i, j \in \nset)$ with $i < j$ give a morphism to the appropriate monoid, and the fact
  that this is an embedding follows from \fullref{Proposition}{prop:staldata} and observing from \eqref{eq:sij-cases}
  that in the image of an element $w$, for $i < j$ (respectively $j < i$) the block corresponding to $s_{ij}$ is equal
  to $I$ (respectively $J$) if and only only if $w$ contains $i$ and factors as $w = w'iw''$ with $j \notin \supp{w''}$.
\end{proof}

\begin{theorem}
  \label{thm_stal_tropical}
  Let $S$ be a commutative unital semiring with zero containing an element of infinite multiplicative order.  The
  stalactic monoid of rank $n$ admits a faithful representation by upper triangular matrices of size $n^2$ over $S$
  having block-diagonal structure with largest block of size $2$ (or size $1$ if $n=1$).
\end{theorem}

\begin{proof}
  Since $(\nset, +)$ embeds in $\ut{1}{S}$ and $\mathcal{F}$ by definition embeds in $\ut{2}{S}$, this is immediate from
  \fullref{Theorem}{thm_stal}.
\end{proof}

\subsection{Variety generated by the stalactic monoid}
\label{subsec:stal-identities}
The finite basis problem for the stalactic monoid has been studied by Han and Zhang \cite[Theorem~4.2]{Han}, but we deduce here some more information about the corresponding variety.

For each word $u$ over a finite alphabet, define $\sigma_u: \supp{u} \rightarrow [|\supp{u}|]$ to be the bijection
taking each symbol $x$ to the number of distinct symbols appearing in the shortest suffix of $u$ containing $x$, that
is, the position of $x$ in an ordering of $\supp{u}$ according to first occurrence when reading $u$ from right to-left.
For example, if $u=dfeebdbf$ is a word over alphabet $\set{a,b,c,d,e,f}$, then $\sigma_u = \begin{psmallmatrix}
    b & d & e & f \\
    2 & 3 & 4 & 1
\end{psmallmatrix}$.

\begin{remark}
  \label{rem:varfgenbyflipflop}
  Let $F$ denote the set of identities $u=v$ with the property that $\sigma_u=\sigma_v$ (and so in particular, the two
  words must have the same support). It is straightforward to verify that the identities $x^2=x, xyx=yx$ form a basis
  for those in $F$ and hence $F$ defines the variety of right regular bands $\var{RRB}$ (see for example
  \cite{petrich}).  Moreover, the variety $\var{RRB}$ is generated by the flip-flop monoid $\mathcal{F}$ (see
  \cite[Proposition 7.3.2]{rhodesteinberg}).
\end{remark}

Recall that $\var{Comm}$ denotes the variety of commutative monoids, and $\var{B}$ the variety of monoids generated by
the bicyclic monoid.

\begin{corollary}
  \label{corol:stalvar}
  Let $n\geq 2$ be a fixed positive integer. The variety of monoids generated by the stalactic monoid of rank $n$ is:
  \begin{enumerate}
  \item a proper subvariety of $\var{B}$;
  \item the join of $\var{Comm}$ and $\var{RRB}$;
  \item defined by the single identity $xyx=yxx$;
  \item equal to the variety generated by the (infinite-rank) stalactic monoid.
  \end{enumerate}
\end{corollary}

\begin{proof}
  \begin{enumerate}
  \item By \fullref{Theorem}{thm_stal}, the stalactic monoid of rank $n$ embeds in the direct product of copies of
    $(\nset_0, +)$ (which embeds in $\uttrop{1}$) and copies of $\uttrop{2}$. Thus $\stal_n$ is contained in the variety
    generated by $\uttrop{2}$. By \cite{daviaud_identities}, the latter is equal to the variety generated by the
    bicyclic monoid. The fact that $\stal_n$ generates a proper variety of the bicyclic variety now follows from the
    fact that Adjan's identity is a minimal length identity for $\mathcal{B}$ \cite{adjan}, whilst $\stal_n$ satisfies
    the identity $xyx=yxx$ \cite[Proposition 15]{cm_identities}.

  \item We first show that the identities satisfied by the stalactic monoid of rank $n$ are precisely the balanced
    identities satisfied by the flip-flop monoid $\mathcal{F}$. By \fullref{Theorem}{thm_stal}, it is clear that the
    stalactic monoid of rank $n$ embeds in the direct product of copies of $(\nset_0, +)$ and copies of
    $\mathcal{F}$. Thus $\stal_n$ satisfies every identity satisfied by both $(\nset_0, +)$ and $\mathcal{F}$. The
    identities satisfied by $(\nset_0, +)$ are precisely identities of the form $u=v$ where $w$ and $v$ have the same
    content, thus $\stal_n$ satisfies all balanced identities satisfied by $\mathcal{F}$. On the other hand,
    $\mathcal{F}$ is an image of $\stal_n$ under any of the homomorphisms $s_{ij}$, and so $\mathcal{F}$ satisfies every
    identity satsified by $\stal_n$.  Finally, note that the identities satisfied by $\stal_n$ are balanced. It thus
    follows that the identities satisfied by the stalactic monoid of rank $n$ are precisely the balanced identities
    satisfied by the flip-flop monoid $\mathcal{F}$. Hence the variety generated by $\stal_n$ is the join of
    $\var{Comm}$ and the variety generated by $\mathcal{F}$. The result now follows from
    \fullref{Remark}{rem:varfgenbyflipflop}.

  \item This has been established in \cite[Theorem~4.2]{Han}, but we give an alternative proof here. We show that each balanced identity in $F$ can be deduced from $xyx=yxx$. Let $u=v$ be such an identity, and
    suppose that $x$ is the rightmost letter of $u$. By repeatedly applying the identity $xyx=yxx$ where $y$ is a factor
    of $u$ lying between two symbols $x$, one can move all symbols $x$ to the right. Inductively it follows that for any
    word $u$ we can deduce from the identity $xyx=yxx$ an identity of the form
    $u =(\sigma_u^{-1}(\ell))^{\alpha_\ell} \cdots (\sigma_u^{-1}(1))^{\alpha_1}$, where $\ell = |\supp{u}|$ and the
    exponents $\alpha_i$ are determined by the content of $u$. Since $\sigma_u=\sigma_v$ and $u$ and $v$ have the same
    content, we may therefore also deduce $u=v$. This shows that each stalactic monoid of rank at least $2$ generates
    the same variety, namely the variety defined by $xyx=yxx$.

  \item Finally, it is clear that $\stal_n$ is contained in the variety generated by $\stal$. Since $\stal$ satisfies the identity $xyx=yxx$, this completes the proof. \qedhere
  \end{enumerate}
\end{proof}

\section{Taiga monoid}
\label{sec:taig}

\subsection{The taiga monoid}

A \defterm{binary search tree with multiplicities} is an ordered (children of each vertex being designated \textit{left}
and \textit{right}), rooted (unless empty) binary tree in which:
\begin{itemize}
\item each vertex is labelled by a positive integer
\item distinct vertices have distinct labels,
\item the label of each vertex is greater than the label of every vertex in its left subtree, and less than the label of
  every vertex in its right subtree
\item each vertex label is assigned a positive integer called its \defterm{multiplicity}.
\end{itemize}
An example of a binary search tree with multiplicities is:
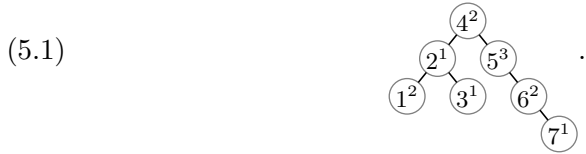
\begin{equation}
\label{eq:bstmulteg}
\begin{tikzpicture}[tinybst,baseline=(0)]
  \node (root) {$4^2$}
    child { node (0) {$2^1$}
      child { node (00) {$1^2$} }
      child { node (01) {$3^1$} }
    }
    child { node (1) {$5^3$}
      child[missing]
      child { node (11) {$6^2$}
        child[missing]
        child { node (110) {$7^1$} }
      }
    };
\end{tikzpicture}.
\end{equation}
The superscript on the label in a vertex denotes its multiplicity.

\begin{algorithm}
\strut\par\nobreak
\textit{Input:} A binary search tree with multiplicities $T$ and a symbol $a \in \nset$.

\textit{Output:} A binary search tree with multiplicities $T \leftarrow a$.

\textit{Method:} If $T$ is empty, create a vertex, label it by $a$, and assign it multiplicity $1$. If $T$ is non-empty,
examine the label $x$ of the root vertex; if $a < x$, recursively insert $a$ into the left subtree of the root; if $a
> x$, recursively insert $a$ into the right subtree of the root; if $a=x$, increment by $1$ the multiplicity of the vertex
label $x$.
\end{algorithm}

Thus one can compute, for any word $u \in \nset^*$, a binary search tree with multiplicities $\ptaig{u}$ by starting with an
empty binary search tree with multiplicities and successively inserting the symbols of $u$, \emph{proceeding right-to-left}
through the word. For example $\ptaig{65117563254}$ is \eqref{eq:bstmulteg}.

Define the relation $\taigcong$ by
\[
u \taigcong v \iff \ptaig{u} = \ptaig{v},
\]
for all $u,v \in \nset^*$. The relation $\taigcong$ is a congruence, and the \defterm{taiga monoid}, denoted $\taig$, is
the factor monoid $\nset^*\!/{\taigcong}$; the \defterm{taiga monoid of rank $n$}, denoted $\taig_n$, is the factor monoid
$[n]^*/{\taigcong}$ (with the natural restriction of $\taigcong$). Each element $[u]_{\taigcong}$ can be identified
with the binary search tree with multiplicities $\ptaig{u}$.

\begin{remark}
  \label{rem:taignodes}
  For a word $u \in \nset^*$ and a symbol $k \in \supp{u}$ there is a unique node of $\ptaig{u}$ containing the symbol
  $k$ with multiplicity $|u|_k$. The number of vertices of $\ptaig{u}$ is therefore equal to $|\supp{u}|$. In a similar
  way to the stalactic monoid (where the structure of $\pstal{u}$ is determined by the content of the word and the order
  of the rightmost instances of the symbols that appear in $u$), it is clear that the structure of the tree $\ptaig{u}$
  is determined by the content of the word and the order of the rightmost instances of the symbols that appear in $u$. The
  difference here is that the construction of the tree also takes into account the ordering of these symbols in the
  underlying alphabet $\nset$.
\end{remark}

\subsection{Constructing a faithful representation of the taiga monoid}

Recall from \fullref{Subsection}{subsec:stal-identities} that for each $u \in \nset^*$ the bijection $\sigma_u : \supp{u} \to [|\supp{u}|]$
is defined to take each symbol to the number of distinct symbols appearing in the shortest suffix of $u$ containing that letter.

\begin{lemma}
\label{taig_subtree}
Let $u \in \nset^*$ and $i,j \in \supp{u}$ with $i \neq j$. The following are equivalent:
\begin{enumerate}
\item $i$ occurs in a subtree of $j$ in $\ptaig{u}$;
\item $\sigma_u(j) < \sigma_u(i)$ and there does not exist $k \in [[i,j]]$ with $\sigma_u(k) <\sigma_u(j)$;
\item $u$ factorizes as $u=u'iu''ju'''$ where $u',u'',u''' \in \nset^*$, $\supp{u'''} \cap [[i,j]] = \emptyset$.
\end{enumerate}
\end{lemma}

\begin{proof}
  By the insertion algorithm, $i$ is in a subtree of $j$ if and only if $i$ is inserted after $j$ and no symbol between
  $k$ with $i < k < j$ was inserted before $j$. The statement (1) is the left-hand side of this equivalence; (2) and (3)
  are different formulations of the right-hand side.
\end{proof}

For $u \in [n]$ and $1 \leq i \neq j \leq n$, let $T_{i,j}(u)$ denote the statement `$u$ factorizes as $u=u'iu''ju'''$
where $u',u'',u''' \in \nset^*$, $\supp{u'''} \cap [[i,j]] = \emptyset$'. Notice that if $i$ and $j$ are in the support
of $u$, then at most one of $T_{i,j}(u)$ and $T_{j,i}(u)$ is true. Otherwise both statements are false.

\begin{proposition}
  \label{taig_data}
  Let $n$ be a fixed positive integer and $u,v \in [n]^*$. The binary search trees with multiplicities $\ptaig{u}$ and
  $\ptaig{v}$ are equal if and only if:
  \begin{enumerate}
  \item $u$ and $v$ have the same content; and
  \item for $ 1 \leq i \neq j \leq n$, $T_{i,j}(u) \iff T_{i,j}(v)$.
  \end{enumerate}
\end{proposition}

\begin{proof}
  If $\ptaig{u}=\ptaig{v}$, then in particular both trees contain the same number of each symbol (hence (1) holds) and
  both trees have the same parent-child structure (and hence by \fullref{Lemma}{taig_subtree} (2) also holds).

  Conversely, suppose first that conditions (1) and (2) hold. By \fullref{Lemma}{taig_subtree} we have that $i$ occurs
  in a subtree of $j$ in $\ptaig{u}$ if and only if $i$ occurs as a subtree of $j$ in $\ptaig{v}$. (Note, if $i<j$, it
  will be the left subtree; if $i>j$, it will be the right subtree.) Thus the statements in (2) determine the same
  parent-child structure of the two trees (the ordering of the children of a given vertex is determined by the order on
  $\supp{u} \subseteq [n]$), while the content of the two words determines the multiplicities of the symbols in each
  trees. Thus (1) and (2) together imply that $\ptaig{u}=\ptaig{v}$.
\end{proof}

Two words in the same taiga class of rank $n$ will have the same image under the content map $c_n$ from
\fullref{Subsection}{subsec:trop-rep}. Therefore, by a slight abuse of notation, consider $c_n$ as a homomorphism from
$\taig_n$ to $\ut{n}{S}$.

Let
\begin{alignat*}{2}
  I &=
  \begin{bmatrix}
    1_S & 1_S & 0_S \\
    0_S & 0_S & 0_S\\
    0_S & 0_S & 0_S
  \end{bmatrix}, &\quad
  J &=
  \begin{bmatrix}
    1_S & 0_S & 0_S\\
    0_S & 1_S & 1_S\\
    0_S & 0_S & 0_S
  \end{bmatrix}, \displaybreak[0]\\
  K &=
  \begin{bmatrix}
    1_S & 0_S & 0_S\\
    0_S & 1_S & 0_S\\
    0_S & 0_S & 0_S
  \end{bmatrix}, &
 L &=
   \begin{bmatrix}
    1_S & 1_S & 1_S\\
    0_S & 0_S & 0_S\\
    0_S & 0_S & 0_S
  \end{bmatrix},
\end{alignat*}
and let $E$ be the $3\times3$ identity matrix. Straightforward calculation shows that $\mathcal{T} = \set{E,K,J,I,L}$ is
a submonoid of $\ut{3}{S}$ with multiplication table as given in \fullref{Table}{tbl:semigroupT}.

\begin{table}[t]
  \caption{The multiplication table of the monoid $\mathcal{T}$.}
  \label{tbl:semigroupT}
  \begin{tabular}{r|ccccc}
        & $E$ & $K$ & $J$ & $I$ & $L$ \\
    \hline
    $E$ & $E$ & $K$ & $J$ & $I$ & $L$ \\
    $K$ & $K$ & $K$ & $J$ & $I$ & $L$ \\
    $J$ & $J$ & $K$ & $J$ & $I$ & $L$ \\
    $I$ & $I$ & $I$ & $L$ & $I$ & $L$ \\
    $L$ & $L$ & $I$ & $L$ & $I$ & $L$ \\
  \end{tabular}
\end{table}

For all $i,j \in \nset$ with $i \neq j$, consider the homomorphism $t_{ij}:\nset^* \rightarrow  \mathcal{T}$
defined by $i \mapsto I$, $j \mapsto J$, all letters $k \in [[i,j]]$ map to $K$, and all other letters map to $E$.

\begin{lemma}
  \label{lem:tij-cases}
  Let $w \in \nset^*$. Then
  \begin{enumerate}
  \item $t_{ij}(w)=E$ if and only if $\supp{w} \cap [[i,j]]= \emptyset$;
  \item $t_{ij}(w)=J$ if and only if $w=w'jw''$ where $i \not\in \supp{w'}$ and $\supp{w''} \cap [[i,j]] = \emptyset$;
  \item $t_{ij}(w)=L$ if and only if $w$ can be factored as $w=w'iw''jw'''$ where $\supp{w'''} \cap [[i,j]] = \emptyset$;
  \item $t_{ij}(w)=I$ if and only if $i$ is in the support of $w$ and $w$ cannot be factored as $w= w'iw''jw'''$ where $\supp{w'''} \cap [[i,j]] = \emptyset$; and
  \item $t_{ij}(w)=K$ otherwise.
  \end{enumerate}
\end{lemma}

\begin{proof}
  Let $W = t_{ij}(w)$. It is easy to see that $W_{3,3}=0$ if and only if $\supp{w} \cap [[i,j]] = \emptyset$, and $E$
  is the only element of $\mathcal{T}$ in which this entry is $0$. It is also easy to see that $W_{2,3}=0$ if and only
  if $w=w'jw''$ where $i \not\in \supp{w'}$ and $\supp{w''} \cap [[i,j]] = \emptyset$, and $J$ is the only element of
  $\mathcal{T}$ in which this entry is $0$. A similar computation shows that $W_{1,3}=0$ if and only if $w=w'iw''jw'''$
  where $i \not\in \supp{w''}$ and $\supp{w'''} \cap [[i,j]] = \emptyset$, and $L$ is the only element of $\mathcal{T}$
  in which this entry is $0$. From the multiplication table above it is easy to see that if $i \in \supp{w}$ then either
  $W=L$ or $W=I$. Thus $W=I$ if and only if $i$ is contained in the support of $w$, but $w$ cannot be factored as
  $w= w'iw''jw'''$ where $i \not\in \supp{w''}$ and $\supp{w'''} \cap [[i,j]] = \emptyset$. \end{proof}

\begin{lemma}
  \label{lem:taighoms}
  Let $n\geq 2$ be a fixed positive integer and let $1 \leq i \neq j\leq n$. Each of the maps
  $t_{ij}: \nset^* \rightarrow \mathcal{T}$ factors to give a homomorphism from the taiga monoid of rank $n$ to
  $\mathcal{T}$.
\end{lemma}

\begin{proof}
  Suppose that $u,v \in [n]^*$ are in the same taiga class of rank $n$. Since $u$ and $v$ have the same content,
  by \fullref{Lemma}{lem:tij-cases} either
  both or neither have image $E$. Moreover, it follows from \fullref{Proposition}{taig_data} that $u$ has image $L$ if
  and only if $v$ has images $L$.
  Hence $i \in \supp{w}$ and $t_{ij}(u) \neq L$ if and only if $i \in \supp{v}$ and
  $t_{ij}(v) \neq L$, showing that $u$ has image $I$ if and only if $v$ has image $I$.

  Thus it suffices to show that that $t_{ij}(u)=J$ if and only if $t_{ij}(v)=J$, or equivalently, that either both words
  admit a factorisation of the form in case~(2) of \fullref{Lemma}{lem:tij-cases}, or neither do. To this end note that
  a word admitting a factorisation of this form is equivalent to $j$ being the first symbol inserted from the interval
  $[[i,j]]$ and the symbol $i$ is never inserted. Given two words $u$ and $v$ in the same taiga class, it is clear that
  symbol $i$ is either in the support of both or neither. Suppose then that $i$ is in the support of neither but, with
  the aim of obtaining a contradiction, that $j$ is the first symbol to be inserted from the interval $[[i,j]]$ when
  reading $u$, whilst some symbol $k$ with $i\neq k\neq j$ is the first symbol from $[[i,j]]$ to be inserted when
  reading $v$. Since $u$ and $v$ have the same support, it follows that we may write $u=u'ku''ju'''$ and
  $v=v'jv''kv'''$, where $k \not\in \supp{u''}$, $j \not\in\supp{v''}$, and
  $\supp{u'''} \cap [[k,j]] = \emptyset = \supp{v'''} \cap [[k,j]]$. But then \fullref{Lemma}{taig_subtree} shows that
  in the tree associated to $w$ symbol $k$ occurs in a subtree of symbol $j$, whilst in the tree associated to $v$
  symbol $j$ occurs in a subtree of symbol $k$. This gives the desired contradiction.
\end{proof}

\begin{theorem}
  \label{thm_taig}
  The taiga monoid $\taig_n$ (respectively, $\taig$) embeds into a direct product of $n$ copies (respectively, countably
  infinite copies) of $(\nset_0, +)$ and $n^2-n$ copies (respectively, countably infinite copies) of the finite
  monoid $\mathcal{T}$.
\end{theorem}

\begin{proof}
The direct product of the morphisms $c_i$ for $i \in [n]$ (respectively, $i \in \nset$) and $t_{ij}$ for $i, j \in [n]$ (respectively, $i, j \in \nset)$ with $i \neq j$ gives a morphism
 to the appropriate monoid, and the fact that this is an embedding follows from
  \fullref{Proposition}{taig_data} together with the observation from \fullref{Lemma}{lem:tij-cases} that in the image
  of an element $w$, the block corresponding to the homomorphism $t_{ij}$ equals $L$ if and only if the statement
  $T_{i,j}(u)$ is true.
  \end{proof}

\begin{theorem}
  \label{thm_taig_tropical}
  Let $S$ be a commutative unital semiring with zero containing an element of infinite multiplicative order. The taiga
  monoid of rank $n$ admits a faithful representation by upper triangular matrices of size $3n^2-2n$ over $S$ having
  block-diagonal structure with largest block of size $3$ (or size $1$ if $n=1$).
\end{theorem}

\begin{proof}
  Since $(\nset, +)$ embeds in $\ut{1}{\trop}$ and $\mathcal{T}$ by definition embeds in $\ut{3}{\trop}$, this is
  immediate from \fullref{Theorem}{thm_taig}.
\end{proof}

\subsection{Variety generated by the taiga monoid}

The finite basis problem for the taiga monoid has also been studied by Han and Zhang; in the next corollary we show how our results can be used to give an alternative proof of \cite[Theorem~4.2]{Han}. Recall that $\var{Comm}$ denotes the variety of commutative monoids and that (from
\fullref{Subsection}{subsec:stal-identities}) $\var{RRB}$ denotes the variety right regular bands.

\begin{remark}
  \label{rem:vargenbyfandtequal}
  The variety $\var{RRB}$ is generated by the monoid $\mathcal{T}$.  Straightforward calculation shows that
  $\mathcal{T}$ satisfies the identities $x^2=x$ and $xyx=yx$.
  On the other hand, as can be seen from \fullref{Table}{tbl:semigroupT}, the monoid $\mathcal{T}$ contains the
  submonoid $\set{E,I,L}$, which is a two element right zero semigroup with an identity adjoined and so isomorphic to
  $\mathcal{F}$. Thus $\mathcal{F}$ lies in the variety generated by $\mathcal{T}$. Thus the monoids $\mathcal{F}$ and
  $\mathcal{T}$ generate the same variety $\var{RRB}$.
\end{remark}

\begin{corollary}
  \label{corol:taigvar}
  Let $n \geq 2$ be a fixed positive integer. The monoids $\stal$, $\stal_n$, $\taig_n$, and $\taig$ each generate the
  same variety.
\end{corollary}

\begin{proof}
  First, $\taig$ is a homomorphic image of $\stal$ \cite[\S~5]{priez_lattice} and thus lies in the variety generated by
  $\stal$. Since $\taig_n$ is a submonoid of $\taig$, it also lies in the variety generated by $\stal$.

  On the other hand, the variety generated by $\taig_n$ contains all commutative monoids and the monoid
  $\mathcal{T}$. By \fullref{Remark}{rem:vargenbyfandtequal}, this variety also contains the flip-flop monoid
  $\mathcal{F}$, and so by \fullref{Theorem}{thm_stal} it contains $\stal$.

  Hence $\stal$, $\taig$, and $\taig_n$ all generate the same variety.  Finally, by \fullref{Corollary}{corol:stalvar},
  $\stal_n$ generates the same variety.
\end{proof}

\section{Sylvester and Baxter monoids}
\label{sec:sylv}

\subsection{The sylvester monoid}

A \defterm{right (respectively, left) strict binary search tree} is an ordered (children of each vertex being designated \textit{left} and \textit{right}), rooted (unless empty) binary tree in which
each vertex is
labelled by a positive integer and the label of each vertex is greater than or equal to the label of every vertex in its
left subtree, and strictly less than the label of every vertex in its right subtree (respectively, the label of each
vertex is strictly greater than or equal to the label of every vertex in its left subtree, and less than or equal to the
label of every vertex in its right subtree).

The following are examples of, respectively, left strict and right strict binary search trees:
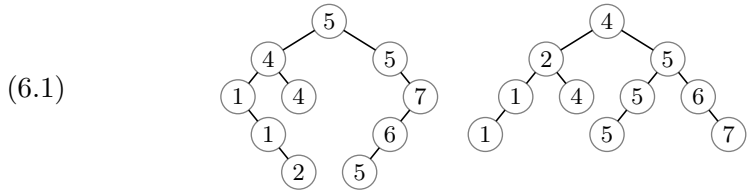
\begin{equation}
  \label{eq:bsteg}
  \begin{tikzpicture}[tinybst,baseline=-10mm]
    \node (root) {$5$}
    child[sibling distance=16mm] { node (0) {$4$}
      child { node (00) {$1$}
        child[missing]
        child { node (001) {$1$}
          child[missing]
          child { node (0011) {$2$} }
        }
      }
      child { node (01) {$4$} }
    }
    child[sibling distance=16mm] { node (1) {$5$}
      child[missing]
      child { node (11) {$7$}
        child { node (110) {$6$}
          child { node (1100) {$5$} }
          child[missing]
        }
        child[missing]
      }
    };
  \end{tikzpicture}
  \quad
  \begin{tikzpicture}[tinybst,baseline=-10mm]
    \node (root) {$4$}
    child[sibling distance=16mm] { node (0) {$2$}
      child { node (00) {$1$}
        child { node (000) {$1$} }
        child[missing]
      }
      child { node (01) {$4$} }
    }
    child[sibling distance=16mm] { node (1) {$5$}
      child { node (10) {$5$}
        child { node (100) {$5$} }
        child[missing]
      }
      child { node (11) {$6$}
        child[missing]
        child { node (111) {$7$} }
      }
    };
  \end{tikzpicture}.
\end{equation}

The insertion algorithms for right (respectively, left) strict binary search trees adds the new symbol as a leaf vertex in the
unique place that maintains the property of being a right (respectively, left) strict binary search tree.

\begin{algorithm}[Right strict leaf insertion]
\label{alg:rightstrictinsertone}
~\par\nobreak
\textit{Input:} A right strict binary search tree $T$ and a symbol $a \in \nset$.

\textit{Output:} A right strict binary search tree $a \rightarrow T$.

\textit{Method:} If $T$ is empty, create a vertex and label it $a$. If $T$ is non-empty, examine the label $x$ of the root; if $a \leq x$, recursively insert $a$ into the left subtree of the root; otherwise recursively insert $a$
into the right subtree of the root. Output the resulting tree.
\end{algorithm}

\begin{algorithm}[Left strict leaf insertion]
\label{alg:leftstrictinsertone}
~\par\nobreak
\textit{Input:} A left strict binary search tree $T$ and a symbol $a \in \nset$.

\textit{Output:} A left strict binary search tree $T \leftarrow a$.

\textit{Method:} If $T$ is empty, create a vertex and label it $a$. If $T$ is non-empty, examine the label $x$ of the root; if $a \geq x$, recursively insert $a$ into the right subtree of the root; otherwise recursively insert $a$
into the left subtree of the root. Output the resulting tree.
\end{algorithm}

Thus one can compute, for any word $u \in \nset^*$, a right strict binary search tree $\psylv{u}$ (respectively, a left
strict binary search tree $\psylvsharp{u}$) by starting with an empty right strict binary search tree (respectively,
left strict binary search tree) and successively inserting the symbols of $u$, \emph{proceeding right-to-left} (respectively,
\emph{left-to-right}) through the word. For example, if $u = 5451761524$, then $\psylvsharp{u}$ and $\psylv{u}$ are
respectively the left- and right-hand trees in \eqref{eq:bsteg}.

Define the relations $\sylvcong$ and $\sylvsharpcong$ by
\begin{align*}
  u \sylvcong v      & \iff \psylv{u} = \psylv{v} \\
  u \sylvsharpcong v & \iff \psylvsharp{u} = \psylvsharp{v},
\end{align*}
for all $u,v \in \nset^*$. The relations $\sylvsharpcong$ and $\sylvcong$ are both congruences.

The \defterm{sylvester monoid}, denoted $\sylv$, is the factor monoid $\nset^*\!/{\sylvcong}$; the \defterm{sylvester
  monoid of rank $n$}, denoted $\sylv_n$, is the factor monoid $[n]^*/{\sylvcong}$ (with the natural restriction of
$\sylvcong$). Each element $[u]_{\sylvcong}$ (where $u \in \nset^*$) can be identified with the binary search tree
$\psylv{u}$.

The \defterm{$\#$-sylvester monoid}, denoted $\sylvsharp$, is the factor monoid $\nset^*\!/{\sylvsharpcong}$; the
\defterm{$\#$-sylvester monoid of rank $n$}, denoted $\sylvsharp_n$, is the factor monoid $[n]^*/{\sylvsharpcong}$ (with
the natural restriction of $\sylvsharpcong$). Each element $[u]_{\sylvsharpcong}$ (where $u \in \nset^*$) can be
identified with the binary search tree $\psylvsharp{u}$.

The monoids $\sylv$ and $\sylvsharp$ are presented (respectively) by $\pres{\nset}{\drel{R}_\sylv}$ and
$\pres{\nset}{\drel{R}_\sylvsharp}$, where
\begin{align*}
  \drel{R}_\sylv &= \gset[\big]{(acvb,cavb)}{a,b,c \in \nset,\; a \leq b < c,\; v \in \nset^*}, \\
  \drel{R}_\sylvsharp &= \gset[\big]{(bvac,bvca)}{a,b,c \in \nset,\; a < b \leq c,\; v \in \nset^*};
\end{align*}
the monoids $\sylv_n$ and $\sylvsharp_n$ are presented (respectively) by $\pres{[n]}{\drel{R}_\sylv}$ and
$\pres{[n]}{\drel{R}_\sylvsharp}$, with the natural restrictions of the sets of defining relations (see \cite[Definition
8]{hivert_sylvester} and \cite[Eq.~(3.10)]{giraudo_baxter2}).

The two monoids $\sylv_n$ and $\sylvsharp_n$ are anti-isomorphic under the map that reflects a binary search tree about a
vertical axis (that is, swaps the left and right children of each node) and replaces each node label $i$ with $n-i+1$.

\subsection{Constructing a faithful representation for the sylvester monoid}
\label{subsec:sylvtroprep}

The first and fourth authors and Ribeiro proved that the sylvester monoid of rank $n$ embeds into $\binom{n}{2}$ copies
of the sylvester monoid of rank $2$, via a map we now define. For $1 \leq i < j \leq n$, define a homomorphism
$\phi_{ij} : [n]^* \to \sylv_2$ by
\[
  a \mapsto
  \begin{cases}
    [1]_{\sylv} & \text{if $a = i$,} \\
    [21]_{\sylv} & \text{if $i < a <j$,} \\
    [2]_{\sylv} & \text{if $a = j$,} \\
    [\emptyword]_{\sylv} & \text{otherwise.} \\
  \end{cases}
\]
Each of these maps factors to give a homomorphism $\phi_{ij} : \sylv_n \to \sylv_2$. Define a homomorphism
$\Phi : \sylv_n \to \prod_{1\leq i < j \leq n} \sylv_2$, where the $(i,j)$-th component of the image of $u \in [n]^*$ is
$\phi_{ij}$.
The map $\Phi$ is injective, as the following result shows:

\begin{proposition} \cite[Lemma 3.6]{cmr_idsylv}
  \label{prop:sylvdata2}
  Let $n$ be a fixed positive integer and $u,v \in [n]^*$. The binary search trees $\psylv{u}$ and $\psylv{v}$ are equal
  if and only if $\phi_{ij}(u) = \phi_{i,j}(v)$ for $i,j \in [n]$ with $i \leq j$).
\end{proposition}

Let
\[
  I = \begin{bmatrix} 1_S&1_S\\ 0_S&0_S\end{bmatrix},\quad
  J = \begin{bmatrix} 1_S&0_S\\ 0_S&\alpha\end{bmatrix},
\]
where $\alpha$ is an element of infinite multiplicative order in the commutative semiring $S$. Then let $\mathcal{M}$ be
the submonoid of $\ut{2}{S}$ generated by $\set{I,J}$. Straightforward calculation shows that $I$ is idempotent,
$JI = I$, and that for any $k \in \nset_0$,
\[
  J^k = \begin{bmatrix} 1_S&0_S\\ 0_S&\alpha^k\end{bmatrix}
  \text{ and }
  IJ^k = \begin{bmatrix} 1_S&\alpha^k\\ 0_S&0_S\end{bmatrix}.
\]
In particular, it is easy to see that $\mathcal{M}$ is the set of all (pairwise distinct) elements $J^k$ and $IJ^k$ for
$k \in \nset_0$, and is presented by $\pres{I, J}{JI=I=I^2}$. From this one can easily observe that $\mathcal{M}$
is isomorphic to a quotient of the sylvester monoid of rank $2$ (since imposing the relations $1^2=1=21$ on
$\sylv_2 = \pres{1, 2}{ 1211=2111, 1221 = 2121}$, yields the monoid $\pres{1, 2}{21=1=1^2}$).

Consider the monoid homomorphism $s : [2]^* \to \ut{2}{S}$ defined by $1 \mapsto I$ and $2 \mapsto J$. Note that the
image of $s$ is $\mathcal{M}$, specifically:
\begin{equation}
  \label{eq:s-cases}
  s(w) = \begin{cases}
    J^k & \text{if $w=2^k$} \\
    IJ^k & \text{if $w=w'12^k$.}
  \end{cases}
\end{equation}

\begin{lemma}
  The homomorphism $s : [2]^* \to \ut{2}{S}$ factors to give a homomorphism from $\sylv_2$ to the monoid $\mathcal{M}$.
\end{lemma}

\begin{proof}
  Let $u,v \in [2]^*$ be such that $u \sylvcong v$. Write $u = u'2^k$ and $v = v'2^j$, where $k$ and $j$ are maximal;
  thus if either $u'$ or $v'$ is non-empty, then it ends with $1$. Then $\psylv{u}$ has exactly $k$ consecutive nodes
  labelled $2$ descending from the root, and $\psylv{v}$ has exactly $j$ consecutive nodes labelled $2$ descending from
  the root. Since $\psylv{u} = \psylv{v}$, it follows that $j = k$ and hence by \eqref{eq:s-cases} we have that
  $s(u) = s(v)$.
\end{proof}

\begin{theorem}
  \label{thm:sylvprod}
  The sylvester monoid $\sylv_n$ (respectively, $\sylv$) embeds into a direct product of $n$ copies (respectively,
  countably infinite copies) of $(\nset_0, +)$ and $\binom{n}{2}$ copies of the (infinite) monoid $\mathcal{M}$.
\end{theorem}

\begin{proof}
  The direct product of morphisms $c_i$ for $i \in [n]$ (respectively, $i \in \nset$) and $s\phi_{i,j}$ for
  $i, j \in [n]$ (respectively, $i,j \in \nset$) with $i<j$ gives a morphism to the appropriate monoid. We show
  that this is an embedding. Let $u, v \in [n]^*$ and suppose that $\psylv{u} \neq \psylv{v}$. By Proposition
  \ref{prop:sylvdata2} there exists $i<j$ such that $\phi_{ij}(u) \neq \phi_{ij}(v)$ in $\sylv_2$.  Since each element
  of $\sylv_{2}$ can be expressed uniquely in the form $2^a1^b2^c$ where $a,b,c \in \nset_0$, we may assume that
  $\phi_{i,j}(u) \sylvcong 2^a1^b2^c \neq 2^p1^q2^r \sylvcong \phi_{i,j}(v)$. If $u$ and $v$ do not have the same
  content, then the fact that $\psylv{u} \neq \psylv{v}$ will be detected in the image of $c_n$. Otherwise, since
  $\phi_{i,j}(w)$ is the element of $\sylv_{2}$ obtained from $w$ by replacing each occurrence of $i$ by $1$; each
  occurrence of $j$ by $2$; each occurrence of a symbol $x$, with $i < x < j$, by $21$; and erasing each occurrence of
  any other element, we see that we must have $a+c=p+r$ and $b=q$ with $c \neq r$. In this case
  $s(\phi_{i,j}(u)) = c \neq r = s\phi_{i,j}(u))$.
\end{proof}

\begin{theorem}
  \label{thm:sylvtroprep}
  Let $S$ be a commutative unital semiring with zero containing an element of infinite multiplicative order. The
  sylvester monoid of rank $n$ admits a faithful representation by upper triangular matrices of size $n^2$ over $S$
  having block-diagonal structure with largest block of size $2$ (or size $1$ if $n=1$).
\end{theorem}

\begin{proof}
  Since $(\nset_0, +)$ embeds in $\ut{1}{S}$ and $\mathcal{M}$ embeds in $\ut{2}{S}$, the result follows by
  \fullref{Theorem}{thm:sylvprod}.
\end{proof}

\begin{remark}
  \label{rem:sylvsharp}
  It is straightforward to prove analogues of \fullref{Theorem}{thm:sylvprod} and \fullref{Theorem}{thm:sylvtroprep} for
  the $\#$-sylvester monoid, using a strategy similar to the above. Specifically, let
  \[ I_{\#} = \begin{bmatrix} \alpha&0_S\\ 0_S&1_S\end{bmatrix},\quad J_{\#} = \begin{bmatrix} 0_S&1_S\\
      0_S&1_S\end{bmatrix}
  \]
  where $\alpha$ is an element of infinite multiplicative order in the commutative semiring $S$. Let $\mathcal{M}_{\#}$
  be the submonoid of $\uttrop{2}$ generated by $\set[\big]{I_{\#},J_{\#}}$. Then one finds that
  \[
    \mathcal{M}_{\#} = \pres[\big]{I_{\#}, J_{\#}}{J_{\#}I_{\#}=J_{\#}=(J_{\#})^2} = \gset{(I_{\#})^k, (I_{\#})^kJ}{k \in \nset_0}
  \]
  is isomorphic to the quotient of $\sylvsharp_2 = \pres{1,2}{2112=2121, 2212=2221}$ obtained by imposing the relations
  $2^2=2=21$. It is clear that $\mathcal{M}_{\#}$ is anti-isomorphic to $\mathcal{M}$.  The monoid homomorphism
  $s_{\#} : [2]^* \to \ut{2}{S}$ defined by $1 \mapsto I_{\#}$ and $2 \mapsto J_{\#}$ has image $\mathcal{M}_{\#}$ with
  \[
    s_{\#}(w) = \begin{cases}
      I_{\#}^k & \text{if $w=1^k$,} \\
      I_{\#}^k J_{\#} & \text{if $w=1^k2w'$.}
    \end{cases}
  \]
  This homomorphism factors to give a homomorphism from $\sylvsharp_2$ to the monoid $\mathcal{M}_{\#}$, and then one
  can construct (i) an embedding of $\sylvsharp_n$ into a direct product of $n$ copies of $(\nset_0,+)$ with
  $\binom{n}{2}$ copies of $\mathcal{M}_{\#}$, and (ii) a faithful upper triangular representation of size $n^2$ for
  $\sylvsharp_n$, in much the same way as in the proofs of \fullref{Theorem}{thm:sylvprod} and
  \fullref{Theorem}{thm:sylvtroprep} (this time making use of the fact that $\sylvsharp_n $ embeds into a direct product
  of copies of $\sylvsharp_2$; see \cite{cmr_idsylv} for details).

  Alternatively, one can see that an analogue of \fullref{Theorem}{thm:sylvtroprep} holds by first observing that the
  permutation $\delta_n: i \mapsto n+1-i$ extends to give an anti-isomorphism
  $\Delta_n: \sylvsharp_n \rightarrow \sylv_n$, and so composing $\Delta_n$ first with the faithful representation given
  by \fullref{Theorem}{thm:sylvtroprep} and then with the anti-isomorphism given by the transpose map will give a
  faithful representation of $\sylvsharp_n$ by \emph{lower} triangular matrices. (From this one can easily obtain an
  upper triangular representation by conjugating the result by the permutation matrix in $\mat{n}{S}$ corresponding to
  the permutation $\delta_n$.)
\end{remark}

\subsection{Variety generated by the sylvester monoid}

The finite basis problem for the sylvester monoid has been studied by Han and Zhang \cite[Theorem~4.6]{Han}; in this section we deduce here some more information about the corresponding variety.

\begin{proposition}
  \label{prop:sylvvar}
  Let $n \geq 2$. Then the variety of monoids generated by the sylvester monoid (respectively, $\#$-sylvester monoid) of
  rank $n$ is:
  \begin{enumerate}
  \item a proper subvariety of $\var{B}$;
  \item the variety generated by the (infinite) monoid $\mathcal{M}$ (respectively, $\mathcal{M}_{\#}$);
  \item not contained in the join of $\var{Comm}$ and any variety generated by a finite monoid;
  \item equal to the variety generated by the (infinite-rank) sylvester monoid (respectively, infinite rank
      $\#$-sylvester monoid).
  \end{enumerate}
\end{proposition}

\begin{proof}
  We prove the result for the case of the sylvester monoids. The arguments follow in a dual manner for the case of the
  the $\#$-sylvester monoids, using \fullref{Remark}{rem:sylvsharp} and the fact that $\sylvsharp_n$ and $\sylv_n$ are
  anti-isomorphic.

  \begin{enumerate}
  \item By \fullref{Theorem}{thm:sylvtroprep}, $\sylv_2$ has a faithful upper-triangular tropical representation with
    maximum block size $2$. Since $\uttrop{2}$ satisfies the same identities as the bicyclic monoid
    \cite{daviaud_identities}, the variety generated by $\sylv_2$ is a subvariety of $\var{B}$. Since $\sylv_n$ embeds
    into a direct product of copies of $\sylv_2$, it follows that the variety generated by $\sylv_n$ is a subvariety of
    $\var{B}$. The containment is proper since $\sylv_n$ also satisfies $xyxy = yxxy$ (see \cite[Proposition
    20]{cm_identities}), while the shortest identity satisfied by the bicyclic monoid has length $10$ (see
    \cite{adjan}).
  \item Since $\mathcal{M}$ is a homomorphic image of $\sylv_2$, which in turn is a homomorphic image of $\sylv_n$, it
    is clear the identities satisfied by $\sylv_n$ must form a subset of the identities satisfied by $\mathcal{M}$. By
    \fullref{Theorem}{thm:sylvprod} the sylvester monoid of rank $n$ embeds in the direct product of copies of
    $(\nset_0, +)$ and copies of $\mathcal{M}$, and so satisfies every balanced identity satisfied by
    $\mathcal{M}$. Since all identities for $\sylv_n$ and $\mathcal{M}$ are balanced, the result follows.
  \item Suppose, with the aim of obtaining a contradiction, that the variety generated by $\sylv_n$ is contained in the
    join of $\var{Comm}$ and a variety generated a finite monoid. Then $\sylv_n$ is an image of a submonoid $L$ of a
    direct product of a commutative monoid and copies of the finite monoid under some surjective homomorphism
    $\phi : L \to \sylv_n$. Let $a,b \in L$ be such that $\phi(a) = [1]_\sylv$ and $\phi(b) = [2]_\sylv$. Then there
    exist $p,q \in \nset$ with $p \neq q$ where $b^p$ and $b^q$ differ only in the commutative monoid component. Thus
    $ab^p$ and $ab^q$ also differ only in the commutative monoid component and so commute with each other. Hence so do
    their images $[12^p]_\sylv$ and $[12^q]_\sylv$. This is a contradiction, since in $\psylv{12^p12^q}$ there are $q$
    nodes $2$ above the topmost node $1$, while in $\psylv{12^q12^p}$ there are $p$ nodes $2$ above the topmost node
    $1$.
  \item This follows immediately from the facts that $\sylv$ embeds into a direct product of copies of $\sylv_2$ and
    that the sylvester monoid of any finite rank is a submonoid of the sylvester monoid of any higher (or infinite)
    rank, as observed in \cite[Theorem 3.9]{cmr_idsylv}. \qedhere
  \end{enumerate}
\end{proof}

While the varieties generated by the hypoplactic, stalactic, and taiga monoids can be expressed as the join of
$\var{Comm}$ and a variety generated by a finite monoid, \fullref[(2)]{Proposition}{prop:sylvvar} shows that this is not
true for a sylvester monoid. Heuristically, this is because there is a bound, dependent only on rank and not on content,
on the number of quasi-ribbon tableaux, stalactic tableaux, and binary search trees with multiplicities with given
content. On the other hand, there is no such bound, independent of content, for binary search trees.

Like the variety generated by the stalatic and taiga monoids, the variety generated by the sylvester monoids turns out to be defined by a single identity, although unlike in those cases the identity in question is not the shortest satisfied by the generating monoid, as the next result (which has been proved by different means in \cite[Theorem
4.14]{cmr_idsylv} and \cite[Theorem
4.6]{Han}) demonstrates.

\begin{theorem}
  \label{thm:sylid}
  The variety of the sylvester monoids is defined by the identity $xyzxty=yxzxty$.
\end{theorem}

\begin{proof}
  We shall refer to the identity $xyzxty=yxzxty$ as the \textit{defining identity} in this proof. By
  \fullref{Proposition}{prop:sylvvar} it suffices to argue using the monoid $\mathcal{M}$.

  First, notice that the defining identity holds in $\mathcal{M}$. Indeed, given any consistent substitution of elements
  of $\mathcal{M}$ for $x$, $y$, $z$ and $t$, it is clear that if $zxty$ maps to an element in the right-zero
  subsemigroup $\gset{IJ^k}{k \geq 0}$ then we have $xyxzty=yxzxty$, otherwise each letter must map to a power of $J$,
  in which case $x$ and $y$ commute.

  Conversely, suppose for a contradiction that some identity holds in $\mathcal{M}$ which is not a consequence of our
  defining identity. Any non-trivial identity can clearly be written uniquely in the form $uaw=vbw$ where $u,v,w$ are
  words and $a \neq b$ are distinct letters.  Choose an identity (over any alphabet), from among those which hold in
  $\mathcal{M}$ but are not consequences of our defining identity, so that when it is written as $uaw=vbw$ as above the
  words $u$ and $v$ are as short as possible. We shall refer to this as the \textit{counterexample identity}.

  First, notice that the letters $a$ and $b$ must both appear in the word $w$, since if not, the substitution sending
  $a$ to $I$, $b$ to $J$ and all other letters to the identity element clearly falsifies our counterexample identity in
  $\mathcal{M}$.

  Since the identity $uaw=vbw$ holds in $\mathcal{M}$, the prefixes $ua$ and $vb$ must have the same content, so in
  particular $u$ contains at least one occurrence of the letter $b$. Write $u=pbq$ with $q$ as short as possible (so in
  particular $q$ does not contain any occurrence of the letter $b$).  Notice that every letter appearing in $q$ also
  appears in $w$, since if not the substitution sending a single letter $c$ of $q$ to $I$, $b$ to $J$ and all other
  letters to the identity element clearly maps $uaw$ to $IJ^{|w|_b}$, whilst $vbw$ is mapped to $IJ^{k}$ for some
  $k>|w|_b$, hence falsifying our counterexample identity in $\mathcal{M}$.

  Now since $a$, $b$ and every letter appearing in $q$ appear in $w$, we may use the defining identity $|qa|$ times
  (with the substitutions mapping $x$ to $b$, $y$ to each letter of $qa$ in turn, $z$ to a suffix of $qa$ followed by a
  prefix of $w$, and $t$ to a subword of $w$) to commute $b$ through $qa$ and deduce the identity $pbqaw = pqabw$. So
  this latter identity is a consequence of the defining identity, and hence also holds in $\mathcal{M}$. Combining with
  the counterexample identity, we see that the identity $pqabw = vbw$ holds in $\mathcal{M}$. The words on either side
  of this identity are the same length as in the counterexample identity and with a longer common suffix, so by the
  minimality assumption on the counterexample identity, the new identity $pqabw = vbw$ must be a consequence of the
  defining identity. But we know that $uaw = pbqaw$ as words, and we have shown that $pbqaw=pqabw$ is a consequence of
  the defining identity, so we deduce that the counterexample identity $uaw=vbw$ is a consequence of the defining
  identity, which is the required contradiction.
\end{proof}

\begin{remark}\label{rem:reverse}
  By \fullref{Proposition}{prop:sylvvar}, the variety of the sylvester (respectively, $\#$-sylvester) monoids is
  generated by $\sylv_2$ (respectively, $\sylvsharp_2$). Since $\sylvsharp_2$ is anti-isomorphic to $\sylv_2$, it is
  clear that reversal of words gives a bijection between the identities satisfied by the sylvester monoid and the
  identities satisfied by the $\#$-sylvester monoids. In particular, the variety of the $\#$-sylvester monoids is
  defined by the single identity $ytxzyx=ytxzxy$, as observed in \cite[Theorem 4.15]{cmr_idsylv}; see also \cite[Theorem 4.7]{Han}.
\end{remark}

\subsection{The Baxter monoid}

For $u \in \nset^*$, define $\pbaxt{u}$ to be the pair $\parens[\big]{\psylvsharp{u},\psylv{u}}$. Define the relation
$\baxtcong$ by
\[
 u \baxtcong v \iff \pbaxt{u} = \pbaxt{v}.
\]
for all $u,v \in \nset^*$. The relation $\baxtcong$ is a congruence. Note that
${\baxtcong} = {\sylvcong} \cap {\sylvsharpcong}$.

The \defterm{Baxter monoid}, denoted $\baxt$, is the factor monoid $\nset^*/\baxtcong$; the \defterm{Baxter monoid of
  rank $n$}, denoted $\baxt_n$, is the factor monoid $[n]^*/\baxtcong$ (with the natural restriction of
$\baxtcong$). The pairs that can arise as $\pbaxt{u}$ for some $u \in \nset^*$ are called pairs of \emph{twin} binary
search trees. This is a proper restriction on the notion of a pair of left- and right-strict binary search trees, but
the details are not important here.

\begin{theorem}
  \label{thm:baxttroprep}
  Let $S$ be a commutative unital semiring with zero containing an element of infinite multiplicative order. The Baxter
  monoid of rank $n$ admits a faithful representation by upper triangular matrices of size $2n^2-n$ over $S$ having
  block-diagonal structure with largest block of size $n$.
\end{theorem}

\begin{proof}
  Both $\sylv_n$ and $\sylvsharp_n$ have faithful representations by upper triangular matrices of size $n^2$,
  each of which is made up of a block of size $n$ (the image of $c_n$) and $\binom{n}{2}$ blocks of size $2$.

  Construct a representation of $\baxt_n$ by taking one copy of the image of $c_n$ and the $2\binom{n}{2}$ blocks of size
  $2$ from the representations of $\sylv_n$ and $\sylvsharp_n$. Since elements of the $\baxt_n$ are equal if and only if
  their projections into $\sylv_n$ and $\sylvsharp_n$ are equal, this representation is faithful.
\end{proof}

\begin{proposition}
  \label{prop:baxtvar}
  Let $n \geq 2$. Then the variety of monoids generated by the Baxter monoid $\baxt_n$ is:
  \begin{enumerate}
  \item the join of the varieties generated by $\sylv_n$ and $\sylvsharp_n$;
  \item the join of the varieties generated by $\mathcal{M}$ and $\mathcal{M}_{\#}$;
  \item a proper subvariety of $\var{B}$;
  \item not contained in the join of $\var{Comm}$ and any variety generated by a finite monoid; and
  \item equal to the variety generated by the (infinite-rank) baxter monoid.
  \end{enumerate}
\end{proposition}

\begin{proof}
	
 \begin{enumerate}
  \item It is easy to see from the definition that an identity is satisfied in $\baxt_n$ if and only if it is satisfied
    in both $\sylv_n$ and $\sylvsharp_n$.
  \item This follows immediately from part (2) of \fullref{Proposition}{prop:sylvvar}.
  \item That the variety generated by $\baxt_n$ is a subvariety of $\var{B}$ follows from part (1) and
    \fullref[(1)]{Proposition}{prop:sylvvar}. To see that it is a proper subvariety, observe that $\baxt$ satisfies the
    identity $xyxyxy = xyyxxy$\cite[Proposition 26]{cm_identities}, while the shortest identity satisfied by semigroups
    in $\var{B}$ has length $10$ on each side (see \cite{adjan}).
  \item This follows from \fullref[(1 and 3)]{Proposition}{prop:sylvvar}.
  \item This follows from the facts that $\baxt$ embeds into a direct product of copies of $\baxt_2$ and that the baxter
    monoid of any finite rank is a submonoid of the Baxter monoid of infinite rank, as shown in \cite[Theorem
    3.12]{cmr_idsylv}. \qedhere
  \end{enumerate}
\end{proof}

Next we provide a finite basis of identities for the variety of Baxter monoids; see also \cite[Theorem 4.16]{cmr_idsylv} or \cite[Theorem~4.10]{Han}
for an alternative approach.

\begin{theorem}
  The variety generated by the Baxter monoids is defined by the two identities:
  \[
    xayb\; xy \; cxdy = xayb\; yx\; cxdy\quad\text{and}\quad xayb\; xy\; cydx = xayb\; yx\; cydx.
  \]
\end{theorem}

\begin{proof}
  We shall refer to the two identities in the statement as the \textit{defining identities}.

  The defining identities can clearly be deduced from the identity $xyzxty = yxzxty$ (which by
  \fullref{Theorem}{thm:sylid} is satisfied in $\sylv$) and separately also from the identity $ytxzyx=ytzxzy$ (which is
  the reverse of the previous identity, and hence is satisfied in $\sylvsharp$ by \fullref{Remark}{rem:reverse}). It
  follows immediately that they are satisfied in $\baxt$.

  Conversely, suppose for a contradiction that some identity is satisfied in $\baxt$ (and hence also in by
  \fullref{Proposition}{prop:baxtvar} in $\mathcal{M}$ and $\mathcal{M}_{\#}$), and is not a consequence of our defining
  identities. As in the proof of \fullref{Theorem}{thm:sylid}, choose such an identity (which will shall call the
  \textit{counterexample identity}) so that when it is written as $uaw=vbw$ with $a \neq b$ single letters the words $u$
  and $v$ are as short as possible, and write $u=pbq$ with $q$ as short as possible. Note that $ua$ and $vb$ must have
  the same evaluation. Since this identity is satisfied in $\mathcal{M}$, we deduce by exactly the same argument as in
  the proof of \fullref{Theorem}{thm:sylid} that $a$, $b$ and every letter of $q$ all occur in $w$.

  We claim that $b$ occurs at least once in $p$. Indeed, if not then since $a \neq b$ and $b$ does not appear in $q$, we
  have that $b$ appears only once in the word $ua$. Since $ua$ has the same evaluation as $vb$, this means that $v$ does
  not contain $b$ and so the morphism mapping $b$ to $J_{\#}$, $a$ to $I_{\#}$ and all other letters to the identity
  element distinguishes $vaw$ and $uab$, hence contradicting the assumption that the counterexample identity holds in
  $\mathcal{M}_{\#}$.

  Next we claim that $a$ and every letter of $q$ also occur in $p$. Indeed, suppose false, say some letter $c$ appears
  in $qa$ but not in $p$. Then the leftmost appearance of $c$ in the word $uaw = pbqaw$ lies within the suffix
  $bqw$. Since the counterexample identity holds in $\mathcal{M}_{\#}$, the leftmost appearance of $c$ in $vbw$ must be
  in the same position (as can be seen by considering the morphism mapping $c$ to $J_{\#}$ and all other letters to
  $I_{\#}$).  Since $ua=pbqa$ and $vb$ have the same content and $q$ does not contain a $b$, this means there are
  strictly more $b$s to the left of the leftmost $c$ in $uaw$ than in $vbw$. But then the morphism mapping $c$ to
  $J_{\#}$, $b$ to $I_{\#}$ and all other elements to the identity distinguishes the words $vaw$ and $uab$, hence
  contradicting the assumption that the counterexample identity holds in $\mathcal{M}$.

  To recap, we have a factorisation (as words) $uaw = pbqaw$ where $a$, $b$ and every letter of $q$ appears in both $p$
  and $w$. Hence we may use the defining identities $|qa|$ times to commute $b$ through $qa$ and deduce the identity
  $pbqaw = pqabw$. So this latter identity is a consequence of the defining identities, and hence also holds in
  $\baxt$. Combining with the counterexample identity, we see that the identity $pqabw = vbw$ holds in $\baxt$. The
  words on either side of this identity are the same length as in the counterexample identity and with a longer common
  suffix, so by the minimality assumption on the counterexample identity, the new identity $pqabw = vbw$ must be a
  consequence of the defining identities. But we know that $uaw = pbqaw$ as words, and we have shown that $pbqaw=pqabw$
  is a consequence of the defining identities, so we deduce that the counterexample identity $uaw=vbw$ is a consequence
  of the defining identities, which is the required contradiction.
\end{proof}

\section{The right patience-sorting monoid}
\label{sec:rps}

\subsection{The right patience-sorting monoid}

An \defterm{rPS-tableau} is a finite array of symbols of $\nset$ in which columns are bottom-aligned, the entries in the bottom row are strictly increasing from left to right, and the entries in each column are non-increasing from top to bottom. For
example,
\[
  \tableau{
    3 \&   \&   \& 6 \\
    2 \&   \& 6 \& 7 \\
    2 \& 4 \& 6 \& 5 \\
    1 \& 3 \& 4 \& 5 \& 7\\
  }
\]
is an rPS-tableau. The insertion algorithm is as follows:

\begin{algorithm}
  \label{alg:rPSinsertone}
  ~\par\nobreak
  \textit{Input:} An rPS-tableau $T$ and a symbol $a \in \nset$.

  \textit{Output:} An rPS-tableau $T\leftarrow a$.

  \textit{Method:} If $a$ is greater than every symbol that appears in the bottom row of $T$, add $a$ to the right of the
  bottom row of $T$. Otherwise, let $C$ be the leftmost column whose bottom-most symbol is greater than or equal to
  $a$. Slide column $C$ up by one space and add $a$ as a new entry of $C$. Output the new tableau.
\end{algorithm}

Thus one can compute, for any word $u \in \nset^*$, an rPS-tableau $\prps{u}$ by starting with an empty rPS-tableau and
successively inserting the symbols of $u$, \emph{proceeding left-to-right} through the word. Define the relation
$\rpscong$ by
\[
u \rpscong v \iff \prps{u} = \prps{v}
\]
for all $u,v \in \nset^*$. The relation $\rpscong$ is a congruence, and the \defterm{right patience sorting monoid},
denoted $\rps$, is the factor monoid $\nset^*\!/{\rpscong}$; the \defterm{right patience sorting monoid of rank $n$},
denoted $\rps_n$, is the factor monoid $[n]^*/{\rpscong}$ (with the natural restriction of $\rpscong$). Each element
$[u]_{\rpscong}$ (where $u \in \nset^*$) can be identified with the rPS-tableau $\prps{u}$.

The possible bottom rows of rPS-tableaux are clearly in bijective correspondence with subsets of $[n]$ (the
correspondence taking each bottom row to the set of elements appearing in it). Let $B$ be the power set of $[n]$, the
elements of which we think of as possible bottom rows of rPS-tableaux, and which we also identify with single-row
rPS-tableaux in the obvious way. (In particular this means we identify the empty set $\emptyset$ with the identity
element in the right patience-sorting monoid.)

\subsection{Constructing a faithful representation}

Note that, unlike the syl\-vester monoid, it is impossible for the right patience sorting monoid of arbitrary rank to be embedded into a direct product of copies of right patience sorting monoid of rank $2$. This is because $\rps_n$ does not
satisfy any identity of length less than $n$ \cite[Proposition~4.8]{cms_patience1}, while it does satisfy the identity
$(xy)^{n+1} = (xy)^nyx$ \cite[Proposition~4.7]{cms_patience1}. Since there is no such embedding, and since rPS-tableaux
are not characterized by content and a bounded amount of extra information, a very different approach is needed to
construct a faithful finite dimensional representation.

If $P$ is an rPS-tableau and $z$ is a generator then, as a direct consequence of the nature of the insertion algorithm,
the bottom row of the rPS-tableau $P \leftarrow z$ is uniquely determined by the combination of the bottom row of $P$
and the generator $z$.  Therefore we may define an action of the free monoid $[n]^*$ on the set $B$, such that for any
rPS-tableau $P$ with bottom row $b \in B$ and any $w \in [n]^*$, the bottom row $b\cdot w$ is the bottom row of the
rPS-tabeau $P \leftarrow w$. Since two words which represent the same element of $\rps_n$ act the same on every tableau,
the action of the free monoid induces an action of the right patience-sorting monoid $\rps_n$, where $b\cdot s$ is
defined to be $b\cdot w$ for some $w \in [n]^*$ with $\prps{w} = s$.

\begin{lemma}
  \label{lem:xyinductive}
  Let $x, y \in [n]$ and $s \in \rps_n$ and let $P$ be an rPS tableau. The number of times generator $x$ occurs in
  column $y$ of $P \leftarrow s$ is uniquely determined by the combination of the element $s$, the bottom row of $P$ and
  the number of times $x$ occurs in column $y$ of $P$.
\end{lemma}

\begin{proof}
  Let $z \in [n]$ and consider the insertion of $z$ into the rPS-tableau $P$. The column into which $z$ gets inserted is
  uniquely determined by $z$ and the bottom row of $P$. Thus, the number of times generator $x$ occurs in column $y$ of
  $P\leftarrow z$ is uniquely determined by the combination of the generator $z$ (in particular whether it equals $x$), the bottom
  row of $P$ and the number of times $x$ occurs in column $y$ of $P$. The bottom row of $P\leftarrow z$ is clearly also uniquely determined
  by the generator $z$ and the bottom row of $P$. The result now follows by induction on the length of a word representing the element $s$.
\end{proof}

We are now ready to define a faithful finite dimensional representation of the right patience-sorting monoid over any commutative unital semiring $S$ containing a zero element and an element $\alpha$ of infinite multiplicative order. For $x,y \in [n]$, define a
function $f_{x,y} : \rps_n \to \matS{B}$ as follows.  For each $s \in \rps_n$ and each $p, q \in B$,
\begin{itemize}
\item if $p\cdot s = q$ then $[f_{x,y}(s)]_{p,q}$ is equal to $\alpha^i$ where $i$ denotes the number of extra $x$s added to column $y$ of any rPS-tableau with bottom row $p$ (for example, $p$ itself viewed as a single-row rPS-tableau) and right-multiplying by $s$.  This is well-defined by \fullref{Lemma}{lem:xyinductive}.
\item if $p\cdot s \neq q$ then $[f_{x,y}(w)]_{p,q} = 0_S$.
\end{itemize}
Notice that, by definition, each row of $f_{x,y}(s)$ contains exactly one non-zero entry; namely the entry in column $p\cdot s$. In particular, the identity element of $\rps_n$ is mapped to the usual identity matrix under this map.
\begin{lemma}
  \label{lemma:rpsmorphism}
Let $S$ be a commutative unital semiring with zero containing an element of infinite multiplicative order. The map $f_{x,y}$ defined above is a morphism from $\rps_n$ to $\matS{B}$.
\end{lemma}

\begin{proof}
  If $u, v \in \rps_n$ then
  \[
    [f_{x,y}(u) f_{x,y}(v)]_{p,q} = \sum_{r \in B} \Big( f_{x,y}(u)_{p,r} \cdot f_{x,y}(v)_{r,q} \Big) = f_{x,y}(u)_{p,p \cdot u} f_{x,y}(v)_{p \cdot u,q}.
  \]
where the left-hand equality is the definition of matrix multiplication, and the right-hand equality is because row $p$ of $f_{x,y}(u)$ contains a non-zero entry only in column $p \cdot u$. Now if $f_{x,y}(v)_{p\cdot u,q} = 0_S$ then this means
  $(p\cdot u)\cdot v = p\cdot (uv) \neq q$, so
  \[
    [f_{x,y}(uv)]_{p,q} = 0_S = [f_{x,y}(u) f_{x,y}(v)]_{p,q}.
  \]
  Otherwise $(p \cdot u) \cdot v = p \cdot (uv) = q$. In
  this case $[f_{x,y}(u)]_{p,p\cdot u}$ is equal to $\alpha^i$ where $i$ denotes the number of $x$s added to column $y$ when taking any rPS-tableau with bottom row
  $p$ and right-multiplying by $u$, and $[f_{x,y}(v)]_{p\cdot u,q} = f_{x,y}(v)_{p\cdot u,p \cdot (uv)}$ is equal to $\alpha^j$ where $j$ denotes the number of $x$s added to column
  $y$ when taking any rPS-tableau with bottom row $p\cdot u$ and right-multiplying by $v$. Thus
    \[
      [f_{x,y}(u) f_{x,y}(v)]_{p,q} = \alpha^i\alpha^j = \alpha^{i+j} = f_{x,y}(uv)_{p,q},
    \]
   since $i+j$ is clearly the number of $x$s added to column $y$ when taking any rPS-tableau with bottom row $p$ and right-multiplying by $uv$, and so $\alpha^{i+j}$ is, by definition, equal to
  $f_{x,y}(uv)_{p,p\cdot (uv)} = f_{x,y}(uv)_{p,q}$.
\end{proof}

\begin{theorem}
  \label{thm:rpstroprep}
  Let $S$ be a commutative unital semiring with zero containing an element of infinite multiplicative order. The right patience sorting monoid of rank $n$ admits a faithful representation by upper-triangular matrices
  of size $ 2^{n-1} (n^2 + n)$ over $S$.
\end{theorem}

\begin{proof}
  Construct a block-diagonal representation by taking the blocks to be the images of the homomorphisms
  $f_{x,y}$ for $x,y \in [n]$ with $x \geq y$.

 To show that the representation is faithful, suppose $u, v \in \rps_n$ are distinct elements of the right
 patience-sorting monoid. Then we may choose $x, y \in [n]$ such that the rPS-tableaux corresponding to these two
 elements differ in the number of times that generator $x$ occurs in column $y$. Clearly we must have $x \geq y$, since
 otherwise neither tableau can contain an $x$ in column $y$. Let $b = \emptyset\cdot u \in B$ be the bottom row of the
 rPS-tableau corresponding to $u$. Then by the definition of $f_{x,y}$ the entry in row $\emptyset$ and column $b$ of
 $f_{x,y}(u)$ contains $\alpha^i$ where $i$ denotes the number of times symbol $x$ occurs in column $y$ of the tableau corresponding to $u$. In
 contrast, the corresponding entry in $f_{x,y}(v)$ is different: either it is $0_S$ (if the tableau of $v$ has a
 different bottom row) or else it is $\alpha^j$ where $j$ denotes the number of times symbol $x$ occurs in column $y$ of
 the tableau corresponding to $v$, where  $i\neq j$ by assumption. Hence, $f_{x,y}(u) \neq f_{x,y}(v)$ so $f(u) \neq f(v)$.

 Notice that the action of $\rps_n$ on the set $B$ has no cycles except for fixed points: right-multiplying an
 rPS-tableau by a generator increases the length of the bottom row, or keeps the length the same and decreases the sum
 along the row, or leaves the row unchanged. It follows that the relation on $B$ defined by $p \leq q$ if and only if
 $p\cdot s = q$ for some $s \in rPS_n$ is a partial order. For each $x$ and $y$ it is immediate from the definition of
 $f_{x,y}$ that for any $w \in [n]^*$ we will have $f_{x,y}(w)_{p,q} = 0_S$ unless $p \leq q$. Thus, completing this
 partial order to a linear order on $B$ yields an order with respect to which $f_{x,y}(w)$ is upper triangular for all
 $w \in [n]^*$. Using such a linear order within each block of the representation thus ensures that the representation is
 upper triangular.

 For each $x$ and $y$, the dimension of the representation $f_{x,y}$ is $|B| = 2^n$.  The block diagonal representation
 has $\binom{n+1}{2}$ blocks of this size and so has total dimension $2^n \binom{n+1}{2} = 2^{n-1} (n^2 + n)$.
\end{proof}

\subsection{Identities satisfied by the right patience-sorting monoid}
An immediate consequence of  \fullref{Theorem}{thm:rpstroprep} is that the right patience sorting monoid $\rps_n$ satisfies every semigroup identity satisfied by $\utS{N}$, where $N=2^{n-1}(n^2+n)$. Taking $S=\trop$ one can therefore determine families of non-trivial identities satisfied by $\rps_n$ by appealing to the known results about $\uttrop{N}$ (see for example \cite{Izhakian}, \cite{Okninski}). In fact, since the representation $f$ constructed in the proof is block-diagonal where each block has size $2^n$, one could instead take $N=2^n$ in the above statement.  Since the variety generated by $\uttrop{N}$ properly contains the variety generated by $\uttrop{d}$ whenever $d < N$ (see \cite[Theorem 2.4]{Aird}), one may be tempted to seek ways to reduce the dimension of the representation,  for the purpose of studying identities.

The dimension of the representation in \fullref{Theorem}{thm:rpstroprep} could be refined down slightly by \textit{ad
  hoc} arguments showing that certain rows and columns are not needed, but any significant reduction in the asymptotics
is likely to require a completely different representation, if it is possible at all. Therefore we instead turn our attention to
the so-called \emph{chain length} of the representation.  Given an $N$-dimensional upper triangular matrix
representation $\phi: M \rightarrow \utS{N}$ of a semigroup $M$ over a semiring $S$, let $\Gamma_\phi$ denote the transitive closure of the directed graph
with node set $\set{1, \ldots, N}$ with an edge from $i$ to $j$ whenever $i=j$ or there exists $A \in \phi(M)$ with $A_{i,j} \neq 0_S$. Writing $i \preceq j$ if there is an edge from $i$ to $j$ in $\Gamma_\phi$, it is clear that $\preceq$ is a partial order. We define the chain length of the representation $\phi$ to be the maximal length of a chain in this partial order. In the case where $S = \trop$ one can use the results of \cite{daviaud_identities} to show that if $\phi: M \rightarrow \uttrop{N}$ is a faihtful tropical representation of chain length $d$, then $M$ is contained in the variety generated by $\uttrop{d}$. To see this, let $\Gamma_\phi(\trop)$ be the set of all matrices $A \in \uttrop{N}$ such that $A_{i,j}=-\infty$ whenever  $i \not\preceq j$ in $\Gamma_\phi$. By construction, the image of the faithful representation $\phi$ is contained in $\Gamma_\phi(\trop)$, and so the statement follows from \cite[Theorem 5.3]{daviaud_identities} which states that the variety generated by $\Gamma_\phi(\trop)$ is equal to the variety generated by $\uttrop{d}$.

We now demonstrate that the chain length of our faithful representation of the right
patience sorting monoid turns out to be far smaller than the dimensions of the blocks.

\begin{proposition}
  The maximum chain length of the representation $f$ is $\binom{n+1}{2}+1$.
\end{proposition}

\begin{proof}
  If $P$ is a tableau and $x$ a generator then the bottom row of $P\leftarrow x$ differs from the bottom row of $P$ if
  and only if insertion of $x$ into $P$ places the new $x$ in a column not already containing an $x$. It follows that
  for $w \in [n]^*$, the number of times the bottom row changes during the iterative construction of the tableau
  $\emptyset \leftarrow w$ is equal to the sum over columns in $\emptyset w$ of the number of distinct generators
  appearing in each column. In particular, this number is an invariant of the element of $\rps_n$ (independent of the
  choice of representative word $w$). Since the $y$th column from the left cannot contain more than $n-y+1$ distinct
  generators (those symbols $x$ with $x \geq y$), this number is bounded above by $\sum_{y=1}^{n} n-y+1 = \binom{n+1}{2}$.

  Now, notice that the block structure of the representation $f$ means that any chain in the partial order $\Gamma_f$ must lie within $\Gamma_{f_{x,y}}$ for some $x$ and $y$ which is clearly the partial order on $B$ given by $p \leq q$ if and only if $p \cdot s = q$ for some $s \in \rps_n$.  Suppose $b_1, \dots, b_k \in B$
  is a chain of distinct sets such that there exist words $w_1, \dots w_{k-1} \in [n]^*$ with
  $f_{x,y}(w_i)_{b_i,b_{i+1}} \neq 0_S$ for each $i \in [k-1]$. By the definition of $f_{x,y}$ this means that
  $b_i \cdot w_i = b_{i+1}$ for each $i \in [k-1]$. Let $w_0 \in [n]^*$ be a word representing the single-row tableau
  $b_1$, and let $w = w_0 w_1 \dots w_{k-1} \in [n]^*$. Consider the corresponding rPS-tableau $\emptyset \leftarrow
  w$. Clearly each $b_i$ appears as a bottom row during the iterative insertion of the symbols in $w$. But the number of
  bottom rows so appearing exceeds by at most $1$ the number of times the bottom row changes. Thus, by the previous
  paragraph. $k \leq \binom{n+1}{2}+1$.

  Conversely, consider the word
  \[
    w = \prod_{i=1}^{n}  \left( \prod_{j=1}^{i} (n-j+1) \right).
  \]
  It is easy to see that the corresponding tableau $\emptyset w$ is of size $\binom{n+1}{2}$ with distinct entries in
  every column. Consider the iterative construction of this tableau using the word $w$. By the first paragraph of the
  proof, the bottom row changes $\binom{n+1}{2}$ times, yielding a sequence of $\binom{n+1}{2}+1$ distinct bottom rows which
  form a chain.
\end{proof}

\begin{corollary}\label{rPStrop}
The right patience-sorting monoid of rank $n$ satisfies all semigroup identities satisfied by $\uttrop{d}$ where $d = \binom{n+1}{2}+1$.
\end{corollary}

\begin{remark} Using some straightforward generalisations of proofs from \cite[Section 5]{daviaud_identities}, one can establish an analogue of \fullref{Corollary}{rPStrop} in which $\trop$ is replaced by any commutative unital semiring $S$ with both a zero element and an element of infinite multiplicative order. However, to our knowledge the only such semirings where much is known about identities in $\utS{d}$ either satisfy no identities (for example, $S = \mathbb{N}$ and hence also $S$ any ring of characteristic $0$) in which case the result is vacuous, or satisfy the same identities as $\trop$ (for example, the tropical natural number semiring $\trop \cap (\mathbb{N} \cup \lbrace -\infty \rbrace$)) in which case the stronger statement adds nothing.
However, we note that the $\rps_n$ certainly satisfies identities which do not come from $UT_d(\trop)$ (for example, those given by \cite[Proposition 4.7]{cms_patience1}, which are shorter than any which hold in $UT_d(\trop)$), so it may be interesting to investigate representations over other semirings.\end{remark}

\begin{remark} The right patience-sorting monoid of rank $n$ does not satisfy any identity of length less than $n$ \cite[Proposition 4.8]{cms_patience1}. Thus, in contrast to the monoids considered in the previous sections,  the identities satisfied by $\rps_n$, and hence the variety generated by $\rps_n$, are dependent on
$n$. Thus there is no direct analogue for right patience-sorting monoids of \fullref{Corollaries}{corol:hypovar}, \ref{corol:stalvar}, or \ref{corol:taigvar} or
		\fullref{Proposition}{prop:sylvvar} (which state that the variety generated by one of our previous plactic-like monoids of infinite rank is equal to the variety generated by any one of the finite rank plactic-like monoids of the same kind). However, there is scope to study other aspects of these varieties further (for example, bases of identities).
\end{remark}

\begin{remark} Since our initial motivation was to study identities for plactic-like monoids via tropical representations, we have not attempted to construct representations of the left-patience sorting monoids; as explained in the introduction, the existence of free submonoids in $\lps_n$ implies that faithful representations by \emph{tropical} matrices simply do not exist. However, one could investigate whether faithful finite-dimensional representations of $\lps_n$ monoids exist over other semirings. A key property of right patience-sorting tableaux used in the construction of our representations is that there are finitely many possible `bottom rows' in an $\rps$-tableau. This is not the case for left-patience sorting monoids, and so a completely different approach to that for $\rps$-monoids would be required.
	\end{remark}

\bibliography{\jobname}
\bibliographystyle{alphaabbrv}

\end{document}